\documentclass[journal]{IEEEtran}
\usepackage{textcomp}
\usepackage{enumitem}
\usepackage{amsthm}
\usepackage{graphicx}
\usepackage{subcaption}
\usepackage{epstopdf}
\usepackage{float}
\usepackage[nocomma]{optidef}
\usepackage{dsfont}
\usepackage{times}
\usepackage{bm}
\usepackage{fancyhdr,graphicx,amsmath,amssymb}
\usepackage[ruled,vlined]{algorithm2e}
\usepackage{hyperref}
\hypersetup{
  colorlinks   = true, 
  urlcolor     = blue, 
  linkcolor    = blue, 
  citecolor    = red   
}
\include{pythonlisting}
\usepackage{algorithmic}
\newcounter{algsubstate}


\def\cu#1{{\color{black}#1}} 
\def\bt#1{{\color{black}#1}}
\def\btt#1{{\color{black}#1}}
\usepackage{cite}
\usepackage{xcolor}
\usepackage{amsmath}
\usepackage{amssymb}
\usepackage{mathtools}

\DeclareMathOperator*{\argmax}{arg\,max}
\makeatletter
\def\BState{\State\hskip-\ALG@thistlm}
\makeatother

\usepackage{stfloats}

\begin{document}

\addtolength{\abovedisplayskip}{-.1cm}
\addtolength{\belowdisplayskip}{-.1cm}

\addtolength{\textfloatsep}{-.5cm}

%
\title{Resilient Primal-Dual Optimization Algorithms for Distributed Resource
Allocation}

\author{\IEEEauthorblockN{Berkay Turan$^\dagger$ \quad}
\IEEEauthorblockN{C\'esar A. Uribe$^\dagger$ \quad}
\IEEEauthorblockN{Hoi-To Wai \quad}
\and
\IEEEauthorblockN{Mahnoosh Alizadeh}
\vspace{-.1cm}}


%


\maketitle


\begin{abstract}
Distributed algorithms for multi-agent resource allocation can provide privacy and scalability over centralized algorithms in many cyber-physical systems. However, the distributed nature of these algorithms can render these systems vulnerable to man-in-the-middle attacks that can lead to non-convergence and infeasibility of resource allocation schemes. In this paper, we propose attack-resilient distributed algorithms based on primal-dual optimization when Byzantine attackers are present in the system. In particular, we design attack-resilient primal-dual algorithms for static and dynamic impersonation attacks by means of robust statistics. For static impersonation attacks, we formulate a robustified optimization model and show that our algorithm guarantees convergence to a neighborhood of the optimal solution of the robustified problem. On the other hand, a robust optimization model is not required for the dynamic impersonation attack scenario and we are able to design an algorithm that is shown to converge to a near-optimal solution of the original problem. We analyze the performances of our algorithms through both theoretical and computational studies.
\end{abstract}

%
\IEEEpeerreviewmaketitle
\newtheorem{proposition}{Proposition}
\newtheorem{corollary}{Corollary}[proposition]
\newtheorem{theorem}{Theorem}
\newtheorem{lemma}{Lemma}
\newtheorem{Fact}{Fact}
\newtheorem{remark}{Remark}
\newtheorem{assumption}{Assumption}
\newtheorem*{runningexample}{Running Example}
\makeatletter
\def\blfootnote{\xdef\@thefnmark{}\@footnotetext}
\makeatother
\newcommand{\prm}{\boldsymbol{\theta}}
\newcommand{\prmdl}{\boldsymbol{z}}
\newcommand{\eqdef}{\vcentcolon=}
\newcommand{\beq}{\begin{equation}}
\newcommand{\eeq}{\end{equation}}
\newcommand{\grd}{\nabla}
\newcommand{\Cset}{\mathcal C}
\newcommand{\ie}{i.e., }
\renewcommand{\thefootnote}{\fnsymbol{footnote}}
 \blfootnote{$^\dagger$Authors have equal contribution. B. Turan and M. Alizadeh are with  Dept. of ECE, UCSB, Santa Barbara, CA, USA. C. A. Uribe is with LIDS, MIT, Cambridge, MA, USA. H. T. Wai is with Dept. of SEEM, CUHK, Shatin, Hong Kong. This work is partially supported by UCOP Grant LFR-18-548175, NSF grant \#1847096, CUHK Direct Grant \#4055113, \cu{and the Yahoo! Research Faculty Engagement Program}. E-mails: \url{bturan@ucsb.edu}, \url{cauribe@mit.edu}, \url{htwai@se.cuhk.edu.hk}, \url{alizadeh@ucsb.edu}}
 \renewcommand{\thefootnote}{\arabic{footnote}}
\section{Introduction}
A number of multi-agent optimization problems arise in a wide range of resource allocation systems that fall under the general umbrella of {\it Network Utility Maximization} problems: in the pioneering example of  congestion control in data networks \cite{kelly,low}; in determining the optimal price of electricity and enabling more efficient demand supply balancing  in smart power distribution systems \cite{mohsenian,lina}; in managing user transmit powers and data rates in wireless cellular networks \cite{mung}; in determining optimal caching policies by content delivery networks \cite{cache}; in optimizing power consumption in wireless sensor networks  with energy-restricted batteries \cite{wsn1,wsn2}; and in designing congestion control systems in urban traffic networks \cite{traffic}.  The shared goal among the above-mentioned problems is to minimize the sum of $N$ user-specific cost functions, subject to a set of  coupling constraints that depend on users' decisions.

In these resource allocation problems, the user-specific cost functions and the set of coupling constraints are considered private information  to the users and to a central coordinator, respectively. Consequently, it is necessary to solve these problems in a distributed fashion \cu{allowing the agents to cooperate through communication with a central coordinator.} Among others, \cu{primal-dual optimization methods}~\cite{koshal2011multiuser} have been advocated as they naturally give rise to decomposable algorithms that favor distributed implementation \cite{palomar2006tutorial}. 
In addition to their practical success, these methods are supported by strong theoretical guarantees where fast convergence to a near-optimal solution is well established \cite{koshal2011multiuser}. 

However, the distributed nature of these methods also exposes the system to vulnerabilities not faced by \cu{their traditional centralized counterpart}. Many of the existing algorithms assume the agents, and the communication channels between the central coordinator and the agents, to be \emph{completely trustworthy}.  
In this paper, we consider the setting where these communications  are susceptible to adversarial attacks. An attacker can take over network sub-systems, and deliberately edit the messages  communicated to the central coordinator to any arbitrary value, \ie a Byzantine attack. As we will demonstrate, this might result in an unstable system with possible damages to hardware and the system overall.

Our goal is to design attack-resilient primal-dual algorithms in order to solve  multi-agent resource allocation problems in presence of Byzantine attackers. If a communication channel is attacked and becomes compromised, the attacker can modify messages and/or inject fresh messages into the network on the agents' behalf. We consider two  scenarios with different attacker capabilities.  A static impersonation attack scenario considers the set of agents communicating through compromised channels to be the same for the duration of the algorithm, whereas a dynamic impersonation attack scenario considers the case where all agents are susceptible to attacks and hence communicate through compromised channels for a {\it limited fraction} of the algorithm's runtime.  
Our main contributions are as follows:
\begin{itemize}[leftmargin=5mm]
    \item We propose resilient distributed resource allocation algorithms under the two aforementioned attack scenarios that rely on robust mean estimation.
    \item \bt{We provide convergence guarantees of the proposed   algorithms.} We show that our algorithm for the dynamic impersonation attack scenario converges to the optimal solution of the regularized problem, while our algorithm for the static impersonation attack scenario converges to an $\mathcal{O}(\alpha_1^2)$ neighborhood of the optimal solution of a robustified and regularized optimization model, where $\alpha_1 \in [0,\frac{1}{2})$ is a known upper bound on fraction of attacked channels.
    \item \cu{We provide empirical evidence that supports our theoretical results on convergence and preventing constraint violation. We do so via computational simulations on electric vehicle charging and power distribution applications.} 
\end{itemize}

\noindent
\textbf{Related work:} Vulnerabilities of various types of distributed algorithms have been identified and addressed in a number of recent studies. Relevant examples \cu{can be found in} \cite{sundaram2011distributed,pasqualetti2012consensus,gentz2016data,sundaram2018distributed,chen2018resilient,consensusbenameur,consensusleblanc,consensusbaras,ravi19consensus} which study secure decentralized algorithms on a  general network topology but consider consensus-based optimization models. There are two fundamental differences between distributed resource allocation and consensus problems that make these algorithms inapplicable in our case:
\begin{itemize}
    \item In resource allocation problems, each agent is solving for their own optimal level of resource consumption, i.e., each agent is solving for their own parameter, whereas consensus problems focus on all agents solving for a shared (global)  parameter.
    \item Unlike resilient consensus algorithms, in resource allocation problems pertaining to access to critical infrastructure systems such as power or transportation networks, one cannot simply block a set of users' access to the network even if they are deemed likely to be attackers.
\end{itemize} 
A recently popular line of works in \cite{feng2014distributed,yin2018byzantine,alistarh2018byzantine,chen2017ml,data19ml} focuses on building resilient algorithms for distributed statistical learning. A crucial difference from this work is that they assume identical functions across the agents.
In fact, we employ robust statistics \cite{Donoho1983,huber2011robust} to develop our resilient algorithms, and particularly, we develop novel results for robust mean estimation, a topic that is recently rekindled in \cite{minsker2015geometric,diakonikolas2016robust,Steinhardt2017}.

The present paper is a revised and extended version of the preliminary conference report \cite{cdcversion}. This paper expands \cite{cdcversion} into multiple attack scenarios and includes numerical studies.

\noindent
\textbf{Paper Organization:}
The remainder of the paper is organized as follows. In Section~\ref{sec:pddra}, we provide an overview of the basic primal-dual algorithm for resource allocation. In Section~\ref{sec:problem}, we formally define two Byzantine attack models and demonstrate how Byzantine attacks can alter the primal-dual optimization procedure. In Section~\ref{sec:algorithms}, we present two attack-resilient primal-dual algorithms corresponding to the different attack scenarios along with their convergence analysis. In Section~\ref{sec:numerical}, we provide numerical results for our algorithms.

\noindent
\textbf{Notations}. Unless otherwise specified, $\| \cdot \|$ denotes the standard Euclidean norm. For any $N \in \mathbb N$, $[N]$ denotes the finite set $\{1,...,N\}$. Given $\prm$, $\prm_i$ indicates the $i$'th block/entry of $\prm$ that corresponds to the parameter of agent $i$. $\prm_{i,j}$ denotes the $j$'th element of vector $\prm_i$. 

\section{Overview of Primal-Dual Algorithm for Resource Allocation}{\label{sec:pddra}}
\algsetup{indent=1em}
\begin{algorithm}[tb]
	\caption{PD-DRA Procedure.}\label{alg:pddra}
	\begin{algorithmic}[1]
        \FOR {$k=1,2,...$}
        \STATE \label{step:accum0} \emph{(Communication stage)}:
        \begin{enumerate}[label=(\alph*)]
        \item \label{step:accum} Central coordinator receives $\{ \prm_i^{(k)} \}_{i=1}^N$ from   agents and computes $\textstyle{\overline{\prm}^{(k)} \eqdef \frac{1}{N} \sum_{i=1}^N \prm_i^{(k)}}$, $\{ \grd_{\prm} g_t (\overline{\prm}^{(k)}) \}_{t=1}^T$.
        \item Central coordinator broadcasts the vector $\overline{\bm g}^{(k)} \eqdef \sum_{t=1}^T \lambda_t^{(k)} \grd_{\prm} g_t (\overline{\prm}^{(k)}) $ to   agents.
        \end{enumerate}
        \STATE \label{step:cn} \emph{(Computation stage)}:
        \begin{enumerate}[label=(\alph*)]
        \item \label{step:parallel} Agent $i$ computes the update for $\prm_i^{(k+1)}$ according to \eqref{eq:pd_prm} using the received $\overline{\bm g}^{(k)} $.
\item The central coordinator computes the update for $\bm{\lambda}^{(k+1)}$ according to \eqref{eq:pd_lam}. 
        \end{enumerate}
        \ENDFOR
	\end{algorithmic} 
\end{algorithm}

We consider the following multi-agent optimization problem with an objective to minimize the average cost incurred by the agents, subject to a set of constraints that are functions of the average of the agents' parameters:
\beq \label{eq:mainopt}
\begin{aligned}
&\underset{ \prm_i \in \mathbb R^d,  \forall i}{\min} & &\bt{f}(\prm)\eqdef\frac{1}{N} \sum_{i=1}^N \bt{f_i} ( \prm_i )\\
&\textnormal{subject to} & & g_t \left( {\textstyle \frac{1}{N} \sum_{i=1}^N \prm_i } \right) \leq 0, &t=1,\dots,T,\\
&&&\prm_i\in\Cset_i, &i=1,\dots,N,
\end{aligned}
\eeq
where $\bt{f_i}(\cdot):\mathbb R^d\rightarrow\mathbb R$ is the continuously differentiable and convex cost function of agent $i$ and $g_t(\cdot):\mathbb R^d\rightarrow\mathbb R$ are continuously differentiable and convex set of constraints. The parameter $\prm_i$ of agent $i$ is constrained to be in a compact convex set $\Cset_i\in\mathbb R^d$. 
\bt{\begin{runningexample}[Resource Allocation Problem]
Throughout the paper, we use the following toy example as a running example to clarify the concepts and the methods: We consider an EV charging example with 5 agents. The cost function $f_i(\cdot)$ is monotone decreasing and is the same for all agents. \btt{As an example, we set $f_i(\prm_i)=(\prm_i-10)^2$ as the quadratic cost function which is monotonically decreasing for $0 \leq \prm_i\leq 10$.} There is a charging station with 5 EV charging points, three of which have a maximum charging rate of $7$kW, and two have a rate of $10$kW. The total rate at which the charging station is able to deliver electricity is determined by the grid, and let it be upper bounded by $25$kW (hence, the average rate is upper bounded by $\frac{25}{5}=5$kW). Accordingly, the constraints of this system are stated as:
\begin{align}
	&\notag \textstyle g\big( (1/5)\sum_{i=1}^5 \prm_i \big)\eqdef (1/5)\sum_{i=1}^5 \prm_i - 5 \leq 0,\\
	&\notag 0\leq\prm_i\leq 7,\quad i=1,2,3,\\
	&\notag 0\leq\prm_i\leq 10,\quad i=4,5.
\end{align}
Note that $\prm$ is a real number, hence dimension $d=1$. The optimal solution in this example is to deliver electricity at a rate of $5$kW to all agents due to symmetry.
\end{runningexample}}
The optimization problem in~\eqref{eq:mainopt} can not be solved centrally, because the utility functions $\bt{f_i}(\cdot)$ are private to the agents, and furthermore the coupling constraints on the resources are only known by a central coordinator. Accordingly, the goal of the primal-dual distributed resource allocation (PD-DRA) procedure in Algorithm~\ref{alg:pddra} is to solve \eqref{eq:mainopt} in a distributed manner, where the agents observe a pricing signal received from the central coordinator and communicate their parameters to the central coordinator \bt{\cite{koshal2011multiuser}}. Consequent to this information exchange, the pricing signal and the agents' parameters are updated by the central coordinator and by the individual agents, respectively. 

In order to derive the update rules used by Algorithm~\ref{alg:pddra}, we first consider the Lagrangian function of \eqref{eq:mainopt}:
\beq
{\cal L} ( \{ \prm_i \}_{i=1}^N; \bm{\lambda} ) \eqdef  
 \frac{1}{N} \sum_{i=1}^N \bt{f_i} ( \prm_i ) + \sum_{t=1}^T \lambda_t \!~ g_t \Big( \frac{1}{N} \sum_{i=1}^N \prm_i \Big),
\eeq
where $\lambda_t\geq 0$ is the dual variable associated with constraint $g_t(\cdot)$ and $\bm \lambda=[\lambda_1\dots\lambda_T]^\intercal\in \mathbb R_+^T$ is the vector of the dual variables. Under strong duality (e.g., when the Slater's condition holds), solving problem~\eqref{eq:mainopt} is equivalent to solving its dual problem:
\beq \label{eq:pd} \tag{P}
\max_{ \bm{\lambda} \in \mathbb R_+^T } ~
\min_{ \prm_i \in \Cset_i, \forall i}~ {\cal L} ( \{ \prm_i \}_{i=1}^N; \bm{\lambda} ) .
\eeq
As suggested in \cite{koshal2011multiuser}, we consider a regularized version of \eqref{eq:pd}. Let us define
\beq \label{eq:reg}
\begin{split}
& {\cal L}_{ \upsilon } ( \{ \prm_i \}_{i=1}^N; \bm{\lambda} ) \eqdef \\
& \textstyle {\cal L} ( \{ \prm_i \}_{i=1}^N; \bm{\lambda} ) + \frac{\upsilon}{2 N} \sum_{i=1}^N \| \prm_i \|^2 - \frac{\upsilon}{2} \| \bm{\lambda} \|^2,
\end{split}
\eeq
such that ${\cal L}_{ \upsilon } (\cdot)$ is $\upsilon$-strongly convex and
$\upsilon$-strongly concave in $\{ \prm_i \}_{i=1}^N$ and $\bm{\lambda}$, respectively.
\bt{
\begin{remark}
\btt{Adding regularization terms is a typical technique used in optimization, called dual smoothing\cite{nesterov2005smooth}. We add the regularization terms for the purposes of convergence analysis used in this paper, which can be applied on strongly convex/concave functions.} Indeed adding the regularization terms might change the solution of the original optimization problem. However, as explained in~\cite[Proposition 5.2]{uribe2018dual}, by an appropriate selection of the regularization parameters, we can recover an optimality gap guarantee for the original problem based on the solution to the regularized problem.
\end{remark}}

We define the regularized problem as: 
\beq \label{eq:pdreg} \tag{P$_\upsilon$}
\max_{ \bm{\lambda} \in \mathbb R_+^T } ~
\min_{ \prm_i \in \Cset_i, \forall i}~ {\cal L_\upsilon} ( \{ \prm_i \}_{i=1}^N; \bm{\lambda} ) .
\eeq

Let $\gamma>0$ be the step size and $k\in \mathbb Z_+$ be the iteration index. The primal-dual recursion performs projected gradient descent/ascent on the primal/dual variables as follows:
\begin{subequations} \label{eq:pda}
\begin{align}
& \prm_i^{(k+1)} = \label{eq:pd_prm} \\
& ~~~~~~{\cal P}_{ \Cset_i } \big( \prm_i^{(k)} - \gamma \!~ \grd_{ \prm_i } 
{\cal L}_{ \upsilon } ( \{ \prm_i^{(k)} \}_{i=1}^N; \bm{\lambda}^{(k)} ) \big),\;\forall i \in [N], \notag \\[.1cm]
& \bm{\lambda}^{(k+1)} = \big[ \bm{\lambda}^{(k)} + \gamma \!~ \grd_{ \bm{\lambda} } {\cal L}_{ \upsilon } ( \{ \prm_i^{(k)} \}_{i=1}^N; \bm{\lambda}^{(k)} ) \big]_+,
\label{eq:pd_lam}
\end{align}
\end{subequations}
where ${\cal P}_{ \Cset_i }(\cdot)$ is the Euclidean projection operator to set $\Cset_i$ and $[\cdot]_+$ denotes $\max\{\cdot,0\}$ operator. According to \eqref{eq:reg}, the gradients with respect to (w.r.t.) the primal and the dual variables are given respectively by:
\begin{subequations}\label{eq:grd}
\begin{equation}\label{eq:grd_prm}
\begin{split}
& \grd_{ \prm_i } {\cal L}_{ \upsilon } ( \{ \prm_i^{(k)} \}_{i=1}^N; \bm{\lambda}^{(k)}   ) = \frac{1}{N} \Big( \grd_{\prm_i} \bt{f_i}( \prm_i^{(k)} ) + \upsilon \!~ \prm_i^{(k)} \\
& \hspace{1.6cm} \textstyle + \sum_{t=1}^T \lambda_t^{(k)} \grd_{\prm} g_t (\prm) \Big|_{\prm = \frac{1}{N} \sum_{i=1}^N \prm_i^{(k)}}  \Big),\\
\end{split}   
\end{equation}
\beq \label{eq:grd_lam}
\begin{split}
& \big[ \grd_{ \bm{\lambda} } {\cal L}_{ \upsilon } ( \{ \prm_i^{(k)} \}_{i=1}^N; \bm{\lambda}^{(k)}  ) \big]_t = g_t \Big( {\textstyle \frac{1}{N} \sum_{i=1}^N} \prm_i^{(k)} \Big) - \upsilon\!~\lambda_t^{(k)},
\end{split}
\eeq
\end{subequations}
for all $i$, $t$. It is worthwhile to highlight that both gradients depend on  the average parameter $\overline{\prm}^{(k)} \eqdef \frac{1}{N} \sum_{i=1}^N \prm_i^{(k)}$. From the above equations \eqref{eq:grd_prm} and \eqref{eq:grd_lam}, we can determine which variables should be communicated between the central coordinator and the agents so that the gradients can be computed locally, see Algorithm~\ref{alg:pddra}. 
\begin{figure*}[t]
\centering
\begin{subfigure}[t]{.8\textwidth}
\centering
   \includegraphics[width=\textwidth]{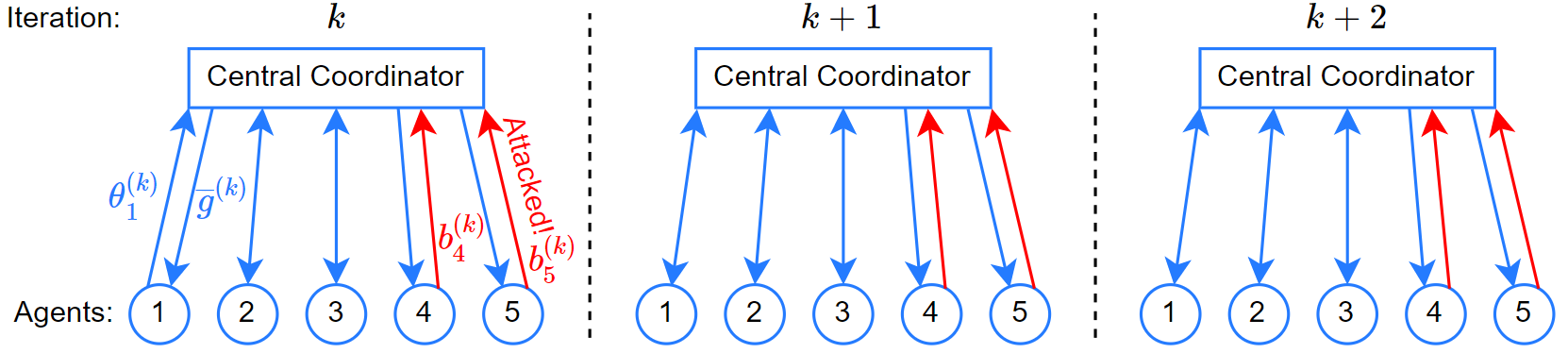}
   \caption{A static impersonation attack scenario, agents 4 and 5 are permanently communicating through compromised channels.}
   \label{fig:permnentattack} 
\end{subfigure}

\begin{subfigure}[t]{.8\textwidth}
   \centering
   \includegraphics[width=\textwidth]{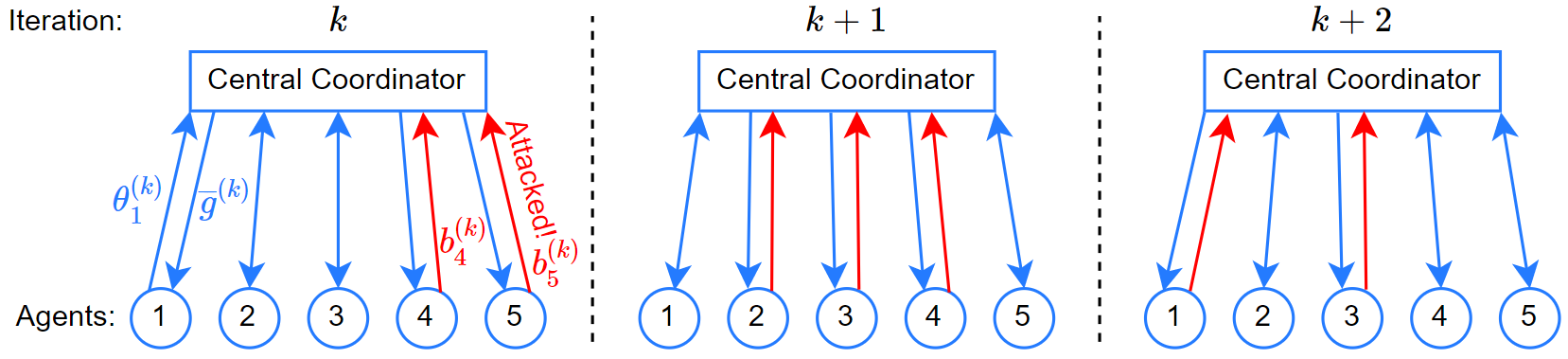}
   \caption{A dynamic impersonation scenario, where the set of agents communicating through compromised channels are changing.}
   \label{fig:Ng2}
\end{subfigure}

\caption{Illustration of (a) static impersonation attack, and (b) dynamic impersonation attack. Blue arrows represent trustworthy channels, whereas red arrows represent compromised channels.}
\label{fig:attackscenarios}

\end{figure*}

Since the regularized primal-dual problem is strongly convex/concave in primal/dual variables, 
Algorithm~\ref{alg:pddra} converges linearly to the optimal solution of \eqref{eq:pdreg} \cite{koshal2011multiuser}. To study
this, let us concatenate the primal and the dual variables and denote ${\bm z}^{(k)} \eqdef ( \{ \prm_i^{(k)} \}_{i=1}^N, \bm{\lambda}^{(k)} )$
as the primal-dual variable at the  $k$th iteration and define the mapping $\bm{\Phi}( {\bm z}^{(k)} )$ as:
\beq\label{eq:mapping}
\bm{\Phi}( {\bm z}^{(k)} ) 
\eqdef \left( \begin{array}{c} 
\grd_{\prm} {\cal L}_{\upsilon} ( \{ \prm_i^{(k)} \}_{i=1}^N, \bm{\lambda}^{(k)}) \\
-\grd_{\bm{\lambda}} {\cal L}_{\upsilon} ( \{ \prm_i^{(k)} \}_{i=1}^N, \bm{\lambda}^{(k)} ) 
\end{array} \right).
\eeq 
\begin{proposition} \cite[Theorem 3.5]{koshal2011multiuser}
Assume that the map $\bm{\Phi}( {\bm z}^{(k)} )$ is $L_{\Phi}$ Lipschitz continuous. 
For all $k \geq 1$, we have
\beq
\| {\bm z}^{(k+1)}  - {\bm z}^\star \|^2 \leq (1 - 2 \gamma \upsilon  + \gamma^2 L_{\Phi}^2 ) \!~ \| {\bm z}^{(k)} - {\bm z}^\star \|^2,
\eeq
where ${\bm z}^\star$ is a saddle point to the \eqref{eq:pdreg}. 
Setting $\gamma = \upsilon  / L_{\Phi}^2$ gives 
$\| {\bm z}^{(k+1)}  - {\bm z}^\star \|^2 \leq \big( 1 - \upsilon^2 / L_{\Phi}^2 \big) \| {\bm z}^{(k)} - {\bm z}^\star \|^2$, $\forall~k\geq 1$.
\end{proposition}

\section{Problem Formulation}\label{sec:problem}

Even though the PD-DRA provides strong theoretical convergence guarantee, it relies on error-free communication between the central coordinator and the agents, and is not robust to attacks on the channels between the agents and the central coordinator, as described below.

We study a situation when the \emph{uplink} communication channels between some of the agents and the central coordinator are  compromised.\footnote{\btt{This paper studies the case where only uplink channels are compromised. However, the case of downlink corruption can also be addressed. Since the downlink channel is a broadcast channel, a compromised downlink channel results in no agent receiving a trustworthy pricing signal. In that case, there is no  optimization method based solution to that problem since there is no communication. If we assume however that all the downlink channels are point-to-point between the central coordinator and each agent, the methods developed in this paper can be applied in a similar fashion.}} Let $\mathcal{A}^{(k)}\subset [N]$ be the set of agents communicating through \emph{compromised uplink channels} at iteration $k$, whose identities are unknown to the central coordinator, and let $\mathcal H^{(k)}\eqdef[N]\setminus {\cal A}^{(k)}$ be the set of agents communicating through \emph{trustworthy uplink channels} at iteration $k$. Instead of receiving $\prm_i^{(k)}$ from each agent $i\in[N]$ at iteration $k$ (Algorithm~\ref{alg:pddra}~Step~\ref{step:accum0}\ref{step:accum}), the central coordinator receives the following messages:
\beq \label{eq:signals}
{\bm r}_i^{(k)} = \begin{cases}
\prm_i^{(k)}, & \text{if}~i \in \mathcal H^{(k)},\\
{\bm b}_i^{(k)}, & \text{if}~i \in \mathcal{A}^{(k)}. 
\end{cases}
\eeq

We consider a Byzantine attack scenario, under which the messages sent through the compromised channels, $\bm b_i^{(k)}$, can be chosen arbitrarily by an adversary. \bt{This also encompasses faulty messages due to erroneous inputs or erroneous channels, since we set no restrictions on $\bm b_i^{(k)}$.} \btt{The adversary's goal is to harm the system and cause suboptimalities. When the messages are erroneous or chosen adversarily, the central coordinator computes the gradients and therefore the pricing signal using these erroneous messages. The agents then update their parameters based on this erroneous pricing signal, which can lead to an overall suboptimal resource allocation.} Moreover, the
choice of the compromised channels ${\cal A}^{(k)}$ affects the impact of the attack and the precautions to be taken in order to defend against the attack. As such, we study two Byzantine attack scenarios that differ in the set of the compromised channels as illustrated in Figure~\ref{fig:attackscenarios}.
\bt{
\begin{runningexample}[Byzantine Attack]
Let agent 1 be communicating through a compromised channel at all iterations, \ie ${\cal A}^{(k)}=\{1\}$, $\forall k$. The compromised message sent to the central coordinator is ${\bm b}^{(k)}_1=1$kW, $\forall k$. This means that irrespective of $\prm_1^{(k)}$, the central coordinator  receives a message indicating agent 1 is willing to charge at rate of $1$kW.
\end{runningexample}
}
\subsection{Attack scenarios}

\begin{enumerate}[wide, labelwidth=!, labelindent=10pt]
\item A \emph{static impersonation attack}, where an adversary takes over a \cu{subset} of uplink channels permanently and the set of agents communicating through compromised channels is fixed (\ie $\mathcal A^{(k)}=\mathcal{A},~\forall k$). Consequently, the central coordinator is never able to communicate reliably with agents $i \in \mathcal A$. In this case, it \cu{is not} feasible to optimize the original problem \eqref{eq:pd} since the contribution from $f(\prm_i): i \in {\cal A}$ becomes unknown to the central coordinator. Yet, we assume that it is also not possible to deny access to resources to agents who are suspected of potentially being under attack. As a compromise, we formulate the following optimization problem:
\beq \label{eq:robust}
\begin{aligned}
&\underset{ \prm_i \in \mathbb \Cset_i, i \in \cal{H}}{\min} & &\bt{f}(\prm)\eqdef\frac{1}{N} \sum_{i\in\cal{H}} \bt{f_i} ( \prm_i )\\
&\textnormal{subject to} & &\underset{\prm_j \in \Cset_j,j\in \mathcal A}{\max} g_t \left( {\textstyle \frac{1}{N} \sum_{i=1}^N \prm_i } \right) \leq 0, &\forall t\in[T].
\end{aligned}
\eeq
The objective of \eqref{eq:robust} is to minimize the cost of the agents with trustworthy channels 
subject to a \emph{robust} set of constraints that consider the \emph{worst case} scenario, \bt{in which the parameters of the agents with compromised channels are assumed to be maximizing the constraints (e.g., those agents are assumed to be consuming the maximum amount of resources). It is critical to mention that during a primal-dual algorithm scheme, the messages received through the compromised channels can be anything. The robust approach is to however ignore those messages, and assume that the parameters of the agents communicating through those channels are maximizing the constraints so that the operation of the system is feasible under any circumstance.} Our goal is to develop an attack-resilient PD-DRA to solve the robust optimization problem \eqref{eq:robust}.
\bt{
\begin{runningexample}[Robust Optimization Model]
Since agent 1 is sending a compromised message of $1$kW and their true parameter can be anything, the worst-case approach is to assume that they are charging at the maximum rate, which is $7$kW for that agent. Hence, the robust constraint is:
\beq
		        \notag
		        \begin{split}
		            &\underset{\prm_j \in \Cset_j,j\in \mathcal A}{\max} g \left( { \frac{1}{5} \sum_{i=1}^5 \prm_i } \right)= \underset{\prm_j \in \Cset_j,j\in \mathcal A}{\max}  \frac{1}{5} \sum_{i=1}^5 \prm_i -5\\
		            &=\frac{1}{5} \sum_{i\in{\cal H}}\prm_i + \underset{\prm_j \in \Cset_j,j\in \mathcal A}{\max}\frac{1}{5}\sum_{j\in {\cal A}}\prm_j -5
		           =\frac{4}{5}\overline{\prm}_{\cal H}-3.6,
		        \end{split}
		        \eeq
		        where we used $|{\cal H}|=4$ and the notation $\overline{\prm}_{\cal H}=\frac{1}{|{\cal H}|}\sum_{i\in{\cal H}}\prm_i$. The robust constraint states that:
\beq\notag
\begin{split}
    \frac{4}{5}\overline{\prm}_{\cal H}-3.6\leq 0\Rightarrow \overline{\prm}_{\cal H}\leq 4.5.
\end{split}
\eeq
The optimal solution in this case is to deliver electricity at a rate of $4.5$kW to the trustworthy agents. Since the compromised agent has the same cost function, their true charging rate will also be $4.5$kW, even though the message sent is $1$kW and the central coordinator assumes their charging rate is $7$kW.
\end{runningexample}
}

\item A \emph{dynamic impersonation attack}, where all the agents \cu{might be} affected by the adversarial attacks but only for a \emph{limited} fraction of time and hence, the set of agents communicating through compromised channels $\mathcal A^{(k)}$ has to dynamically change with iteration $k$. As opposed to the static case, this scenario considers the case where the central coordinator is able to communicate reliably with all the agents at some iterations. Due to this distinction, it is necessary to mention that the static attack is not a special case of the dynamic attack and both scenarios are distinguishable from each other. The dynamic scenario could be applicable when agents do not have dedicated communication channels to the central coordinator and instead communicate over random access systems which are more
appropriate for distributed deployments. Hence, each user periodically accesses authenticated network devices/subsystems that are controlled by Byzantine adversaries and can alter the user's message.  Our goal is to develop an attack-resilient PD-DRA algorithm that can still solve the original regularized problem \eqref{eq:pdreg} in this environment.
\end{enumerate}

\begin{figure}[t]
     \centering
     \hspace{-0.4cm}
     \begin{subfigure}[b]{0.25\textwidth}
         \centering
         \includegraphics[width=\textwidth]{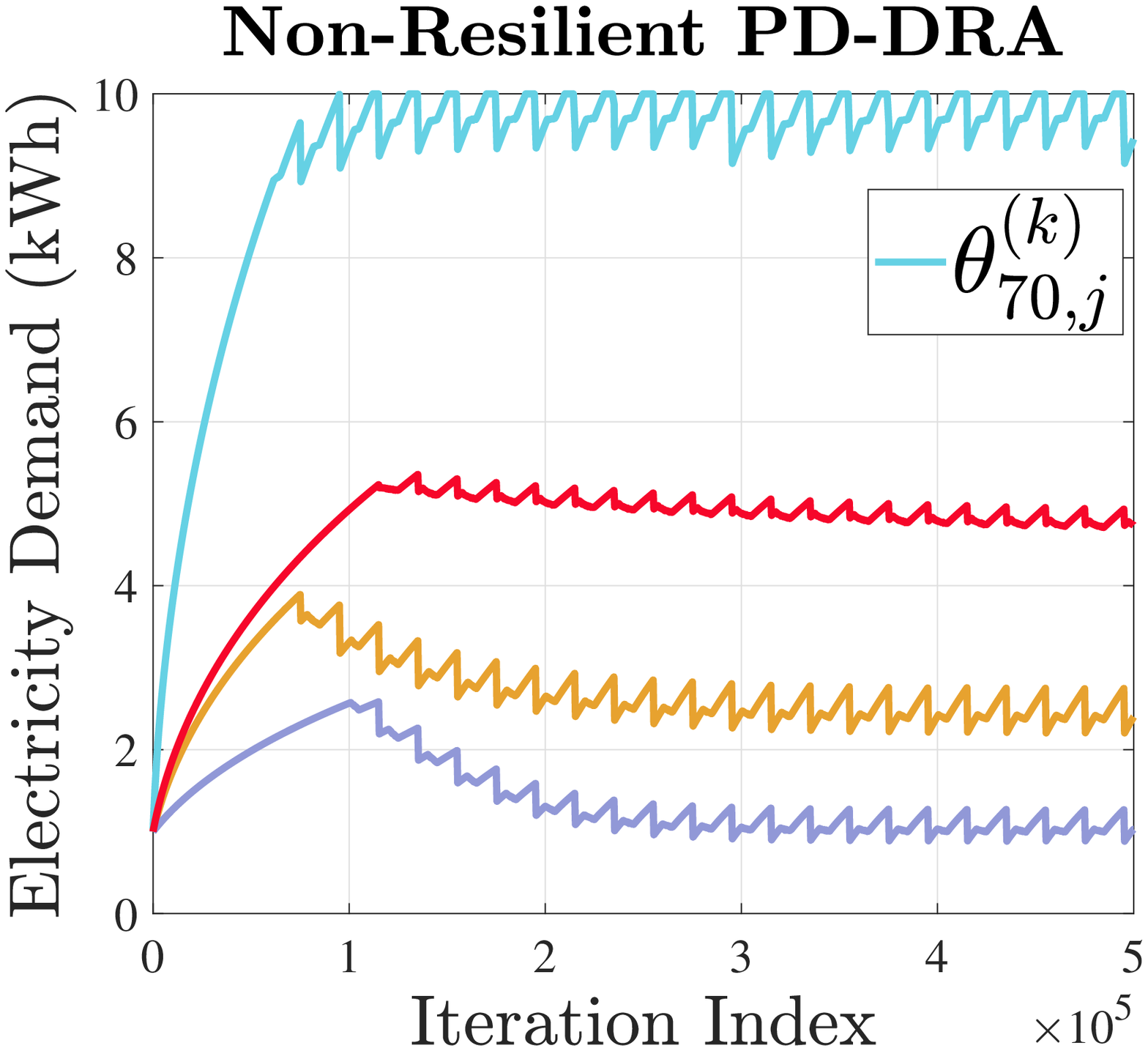}
         \caption{}
         \label{fig:nonresilientparam}
     \end{subfigure}
     \hspace{-0.2cm}
     \begin{subfigure}[b]{0.25\textwidth}
         \centering
         \includegraphics[width=\textwidth]{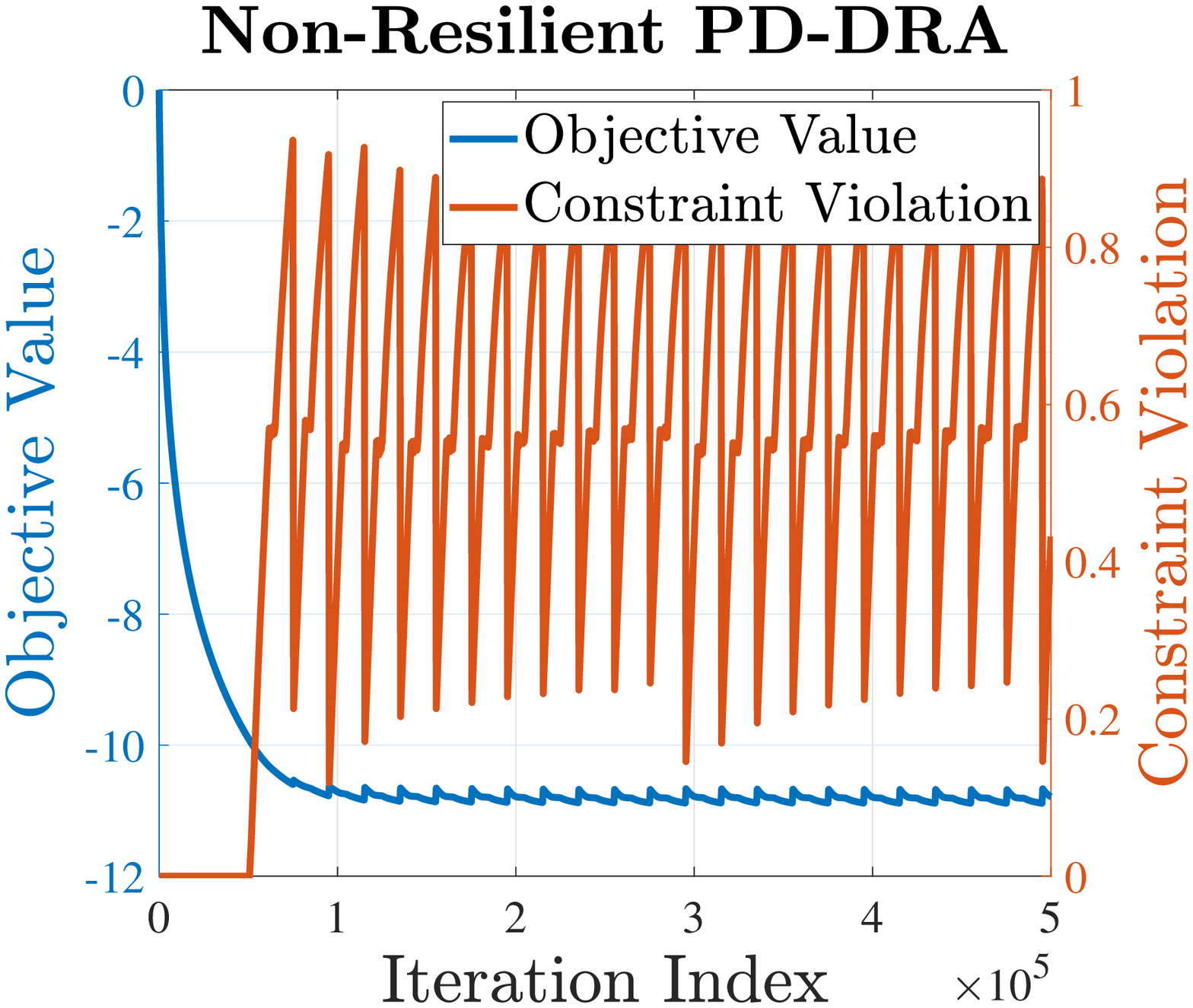}
         \caption{}
         \label{fig:nonresilientobjconst}
     \end{subfigure}
        \caption{Illustration of basic PD-DRA algorithm failure under static impersonation attack. (a) The agents' parameters do not converge, (b) the objective  function does not converge and moreover there is constraint violation. We only display one constraint for brevity.}
        \label{fig:nonresilientfailure}
\end{figure}
\subsection{Limitations of the Basic PD-DRA Algorithm}
Applying the basic PD-DRA algorithm under a Byzantine attack scenario can lead to undesirable outcomes. Recall that the gradients in \eqref{eq:grd} depend on the average parameter $\overline{\prm}^{(k)}$. Under a Byzantine attack scenario, if the central coordinator forms the naive average $\widetilde{\prm}^{(k)} = (1/N) \sum_{i=1}^N {\bm r}_i^{(k)}$ and computes the gradients $\grd g_t( \widetilde{\prm}^{(k)} )$ accordingly, this may result in \btt{large error since the deviation $\widetilde{\prm}^{(k)} - (1/N) \sum_{i=1}^N \prm_i^{(k)}$ can be large (proportional to the maximum diameter of ${\cal C}_i$'s)}. This in turn can obstruct convergence and also overload the system by causing constraint violations. 
\bt{
\begin{runningexample}[Basic PD-DRA Failure]
If the central coordinator believes all the agents are sending trustworthy information, then the optimal solution will occur when one agent is demanding $1$kW and the others are demanding $6$kW (so that the average is $5$kW). But since the $1$kW message is compromised and all the agents have same cost function, the compromised agent's true electricity demand is also at a rate of $6$kW. Hence, the solution delivers electricity at an average rate of 6kW, which is infeasible.
\end{runningexample}
}
We preview our numerical result of applying the basic PD-DRA method under a static impersonation attack scenario for an optimal electric vehicle charging application in Figure~\ref{fig:nonresilientfailure}. For constraint $g_t(\cdot)$, we define constraint violation as $\max\{0,g_t(\overline{\prm}^{(k)})\}$. Observe that the PD-DRA method does not provide convergence and the first constraint is being violated. From resource allocation perspective, this means that the agents are asking to consume more resources than the available amount in the system, which is infeasible. For details regarding the experimental setup, please see Section \ref{sec:numerical}.
\section{Resilient PD-DRA Algorithms}\label{sec:algorithms}
Motivated by the failure of the basic PD-DRA procedure under Byzantine attack scenarios, resilient PD-DRA algorithms are necessary to optimize multi-agent systems in a distributed manner when the system is susceptible to attacks. We hold the following assumption to be true throughout the rest of the paper and  propose two different attack resilient PD-DRA algorithms corresponding to the different attack scenarios outlined in Section~\ref{sec:problem}.
\begin{assumption} \label{ass:bdd}
For all $\prm \in \mathbb R^d$ and for all $t$, the gradient of $g_t$ is bounded with $\| \nabla g_t( \prm ) \| \leq B$ and is $L$-Lipschitz continuous. \bt{Moreover, since maximum resource that can be consumed by an agent is bounded due to limited amount of resources,} we let $\bm 0\in \Cset_i$ and upper bound the diameters of $\Cset_i$ by R:
\beq \label{eq:Rbd}
\max_{ \prm, \prm' \in \Cset_i } \| \prm - \prm' \| \leq R,~i=1,...,N.
\eeq
\end{assumption}
\bt{
\begin{runningexample}[Assumptions]
The constraint in our running example satisfies $\grd g(\prm)=1$, which is bounded by $B=1$ and is $L=0$-Lipschitz continuous. Since the maximum charging rate is upper bounded by $7$kW for three of the agents and by $10$kW for two of the agents, $R=10$. 
\end{runningexample}
}


\subsection{Static Impersonation Attack}\label{sec:static_robust}
Under this attack scenario, given the complete lack of any credible information on the resource consumption parameters of the agents that permanently communicate through compromised channels, the central coordinator can only hope to solve the robust optimization model defined in \eqref{eq:robust} instead. This formulation  considers a worst-case scenario on how much resources the compromised agents will consume, which ensures constraint satisfaction in all cases. However, the constraints in \eqref{eq:robust} require the knowledge of the set $\mathcal A$ 
and the sets $\Cset_j,\forall j\in \mathcal A$, yet the central coordinator lacks this information. 

Hence, in order to develop a robust optimization model that can handle the worst-case scenario without the knowledge of $\cal A$, we let $\alpha_1 \geq |\mathcal A|/N$ as a known upper bound to the fraction of agents communicating through compromised channels and assume $\alpha_1<1/2$, where less than half of the agents are communicating through compromised channels.\footnote{\btt{If more than half of the agents communicate through compromised channels, then the adversary controls the majority and therefore the median, which will be used to estimate the average parameter later in the paper. In that case, there is no optimization based solution the central coordinator can implement in order to securely run the system.}} Let $\overline{\prm}_\mathcal{H}\eqdef\frac{1}{|\mathcal{H}|}\sum_{i\in\mathcal{H}}\prm_i$ be the mean of the agent's parameters that are sent through trustworthy channels. We then define the following set of constraints
\beq
\overline{g}_t ( \prm ) \eqdef g_t( \prm ) + \alpha_1 \big( R B + {\textstyle \frac{1}{2}} L R^2 \big),
\eeq
and formulate a conservative approximation of \eqref{eq:robust}:
\begin{lemma}\label{lem:conservativeform}
Under Assumption \ref{ass:bdd} the following problem yields a conservative approximation of \eqref{eq:robust}, \ie its feasible set is 
a subset of the feasible set of \eqref{eq:robust}:
\beq \label{eq:reformulate_s}
\begin{aligned}
&\underset{ \prm_i \in \mathcal C_i,  i \in \cal{H}}{\min} & &\frac{1}{N} \sum_{i \in \cal{H}} \bt{f_i} ( \prm_i )\\
&\textnormal{subject to} & & \overline{g}_t \left( (1-\alpha_1)\overline{\prm}_\mathcal{H} \right) \leq 0,~\forall~t \in [T],
\end{aligned}
\eeq
\end{lemma}
The proof can be found in Appendix~\ref{app:conservative}.
\bt{\begin{remark}
The proof of Lemma~\ref{lem:conservativeform} is done by upper bounding constraints of \eqref{eq:robust} using Assumption~\ref{ass:bdd} and the fact that $\alpha_1\geq|{\cal A}|/N$. The looser these upper bounds compared to the true values, the more conservative is \eqref{eq:reformulate_s}. This approach potentially leaves less resources available to the agents communicating through trustworthy channels by assuming more than $|{\cal A}|$ number of agents having maximum possible impact on the constraints, irrespective of their set ${\cal C}_i$ or the true value/gradient of the constraints.
\end{remark}}
\bt{
\begin{runningexample}[Conservative Approximation]
With $B=1$, $L=0$, and $R=10$, the conservative approximation of the running example has the following constraint:
\begin{align*}
		        \overline{g}((1-\alpha_1)\overline{\prm}_{\cal H})&=g((1-\alpha_1)\overline{\prm}_{\cal H})+\alpha_1(RB+\frac{1}{2}LR^2)\\
		        &=(1-\alpha_1)\overline{\prm}_{\cal H}-5+10\alpha_1
		    \end{align*}
If $\alpha_1=|{\cal A}|/N=0.2$, then the upper bound is  the fraction of compromised channels. In that case, the constraint is:
\beq\notag
\begin{split}
    0.8\overline{\prm}_{\cal H}-3\leq 0\Leftrightarrow \overline{\prm}_{\cal H}\leq 3.75,
\end{split}
\eeq
which is more conservative compared to the constraint of the robust optimization model (which was $\overline{\prm}_{\cal H}\leq 4.5$). The optimal solution in this case is to deliver electricity at a rate of $3.75$kW to the agents. The conservatism arises due to the difference between agent-specific maximum charging rate $7$kW and the absolute maximum charging rate $10$kW. Since the constraint is linear, the gradient is constant. Hence, the smoothness and Lipschitz bounds hold with equality without causing additional conservatism.\\
If however $\alpha_1=0.4$, then the central coordinator assumes two agents communicating through compromised channels. In this case the conservative approximation has the constraint as
\beq\notag
\begin{split}
    0.6\overline{\prm}_{\cal H}-1\leq 0\Leftrightarrow \overline{\prm}_{\cal H}\leq \frac{5}{3},
\end{split}
\eeq
which results in charging at an even slower rate since the central coordinator has to be robust against two agents charging at the maximum rate of $10$kW.
\end{runningexample}
}
To develop an attack resilient PD-DRA algorithm, we again define the regularized Lagrangian function of $\eqref{eq:reformulate_s}$:
\beq\label{eq:lagrangianrobustreg}
\begin{split}
& \overline{\cal L}_{\upsilon} ( \{ \prm_i \}_{i \in \cal{H}}; \bm{\lambda} ; {\cal H} )  \\
& \eqdef {\textstyle \frac{1}{N} \sum_{i \in \cal{H} }} \bt{f_i} ( \prm_i ) + {\textstyle \sum_{t=1}^T} \lambda_t \!~ \overline{g}_t \left( (1-\alpha_1)\overline{\prm}_\mathcal{H} \right) \\
& \hspace{.5cm} \textstyle + \frac{\upsilon}{2N} \sum_{i \in \cal{H}} \| \prm_i \|^2 - \frac{\upsilon}{2} \| \bm{\lambda} \|^2.
\end{split}
\eeq
The above function is $(1-\alpha_1)\upsilon$-strongly convex and concave in
$\prm$ and $\bm{\lambda}$, respectively 
(since $(1-\alpha_1)\leq\textstyle{\frac{|\mathcal H|}{N}}\leq 1$). 
Our main task is to tackle the following modified problem of \eqref{eq:pd} under Byzantine attack on (some of) the uplinks:
\beq
\label{eq:pdp} \tag{P$'_{\upsilon}$}
\max_{ \bm{\lambda} \in \mathbb R_+^T } 
\min_{ \prm_i \in  \Cset_i, \forall i \in \cal{H}}~ \overline{\cal L}_{\upsilon} ( \{ \prm_i \}_{i \in \cal{H}}; \bm{\lambda}; {\cal H} ) .
\eeq
Notice that \eqref{eq:pdp} bears a similar form as \eqref{eq:pd} and thus  one may  apply the PD-DRA method to the former. The gradients with respect to primal/dual variables are given by:
\begin{subequations}\label{eq:grdrobust}
\begin{equation}\label{eq:grd_prm_robust}
\begin{split}
& \grd_{ \prm_i } \overline{\cal L}_{\upsilon} ( \{ \prm_i^{(k)} \}_{i \in \cal{H}}; \bm{\lambda}^{(k)}; {\cal H} ) = \frac{1}{N} \Big( \grd_{\prm_i} \bt{f_i}( \prm_i^{(k)} ) + \upsilon \!~ \prm_i^{(k)}, \\
&  \textstyle + \frac{(1-\alpha_1)N}{|\mathcal H|}\sum_{t=1}^T \lambda_t^{(k)} \grd_{\prm} \overline{g}_t (\prm) \Big|_{\prm = (1-\alpha_1)\bar{\prm}_\mathcal{H}^{(k)}}  \Big), \btt{\forall i \in {\cal H}},\\
\end{split}   
\end{equation}
\beq \label{eq:grd_lam_robust}
\begin{split}
& \big[ \grd_{ \bm{\lambda} } \overline{\cal L}_{\upsilon} ( \{ \prm_i^{(k)} \}_{i \in \cal{H}};\bm{\lambda}^{(k)}; {\cal H} )\big]_t = \overline{g}_t \Big( (1-\alpha_1)\overline{\prm}^{(k)}_\mathcal{H} \Big) - \upsilon\lambda_t^{(k)}.
\end{split}
\eeq
\end{subequations}
However, such application requires the central coordinator to compute the sample average
\beq \label{eq:truemean}
\textstyle
\overline{\prm}^{(k)}_{\cal H} = \frac{1}{|{\cal H}|} \sum_{ i \in {\cal H} } \prm_i^{(k)},
\eeq 
at each iteration.
The above might not be computationally feasible under the  attack model,
since the central coordinator is oblivious to the identity of ${\cal H}$. 
As a solution, the central coordinator computes the robust mean $\widehat{\prm}_{\mathcal H}^{(k)}$ of the received parameters $\{\bm r_i^{(k)}\}_{i\in[N]}$ using a median-based mean estimator described next.

\subsubsection{Overview of Median-Based Mean Estimation}\label{sec:robustmean}

Consider a set of $N$ vectors $\{{\bm x}_i\in \mathbb{R}^d\}_{i=1}^N$, among which at least $(1-\alpha_1)N$ are trustworthy ($\bm x_i\in \mathcal{H}$) and at most $\alpha_1 N$ are compromised ($\bm x_i\in \mathcal{A}$). We consider a simple median-based estimator applied to each coordinate $j=1,\dots,d$. First, define the coordinate-wise median as:
\begin{equation*}
    \left[ {\bm x}_{\sf med} \right]_j = {\sf med}\left( \{ [{\bm x}_i]_j \}_{i=1}^N \right),
\end{equation*}
where ${\sf med(\cdot)}$ computes the median of the operand. Then, our estimator is computed as the mean of the nearest $(1-\alpha_1)N$ neighbors of $\left[ {\bm x}_{\sf med} \right]_j$. 
Our estimator is:
\begin{equation}\label{eq:median}
    \textstyle
[ \widehat{\bm x}_{\mathcal H} ]_j = \frac{1}{(1-\alpha_1)N} \sum_{i \in {\cal N}_j}  [{\bm x}_i ]_j,
\end{equation}  
where we have defined the set with $|{\cal N}_j| = (1-\alpha_1)N$ as:
\begin{equation*}
    \bt{{\cal N}_j} = \{ i \in [N] : \big| \big[ {\bm x}_i -  {\bm x}_{\sf med} \big]_j \big| \leq r_j \},
\end{equation*}
such that $r_j$ is chosen to satisfy $|{\cal N}_j| = (1-\alpha_1)N$. 

\btt{The following bounds the performance of \eqref{eq:median}:
\begin{proposition}\label{prop:median}
Let $\overline{\bm x}_{\mathcal H}$ be the mean of the trustworthy vectors. Suppose that $\max_{i \in {\cal H}} \| \bm x_i - \overline{\bm x}_{\cal H} \|_\infty \leq r$,
then for any $\alpha_1 \in (0, \frac{1}{2})$, it holds that:
\begin{equation}
    \label{eq:median_bound}
\|\widehat{\bm x}_{\cal H} - \overline{\bm x}_{\cal H} \| \leq \frac{2\alpha_1}{1-\alpha_1} \left( 1 + \sqrt{\frac{(1-\alpha_1)^2}{1-2\alpha_1}} \right) r \sqrt{d}.
\end{equation}
\end{proposition}}


\btt{The proof can be found in Appendix~\ref{app:medianest}. We note that for sufficiently small $\alpha_1$, the right hand side on \eqref{eq:median_bound} can be approximated by ${\cal O}(\alpha_1 r \sqrt{d})$.} Using this median-based mean estimator, we propose the robust PD-DRA algorithm as follows.

\algsetup{indent=1em}
\begin{algorithm}[t]
	\caption{Robust PD-DRA Algorithm}\label{alg:robust_pddra}
	\begin{algorithmic}[1]
		\STATE \textbf{Input}: Each agent has initial state $ \prm_i^{(0)}$.
		\FOR {$k=1,2,...$}
		\STATE \label{r_step:accum0} \emph{(At the Central Coordinator)}:
		\begin{enumerate}[label=(\alph*)]
			\item \label{r_step:accum} Receives $\{{\bm r}_i^{(k)}\}_{i=1}^N$, see~\eqref{eq:signals}, from agents.
			\item \label{r_step:robustmean}Computes robust mean  $\widehat{\prm}_{\cal H}^{(k)}$ using the estimator \eqref{eq:median}.
			\item \label{r_step:broadcast}Broadcasts the vector $\widehat{\bm g}^{(k)}_{\cal H} \eqdef \sum_{t=1}^T \lambda_t^{(k)} \grd_{\prm} \overline{g}_t ((1-\alpha_1)\widehat{\prm}_{\cal H}^{(k)}) $ to   agents.
			\item Computes the update for $\bm{\lambda}^{(k+1)}$  with \eqref{eq:dual_pb}. 
		\end{enumerate}
		\STATE \label{r_step:cn} \emph{(At each agent $i \btt{\in \mathcal{H}}$)}:
		\begin{enumerate}[label=(\alph*)]
			\item Agent receives $\widehat{\bm g}^{(k)}_{\cal H}$.
			\item \label{r_step:parallel} Agent computes update for $\prm_i^{(k+1)}$ with \eqref{eq:primal_pb}.

		\end{enumerate}
		\ENDFOR
	\end{algorithmic} 
\end{algorithm}

\subsubsection{Robust PD-DRA Algorithm}
\begin{figure}[t]
    \centering
    \includegraphics[width=0.4\textwidth]{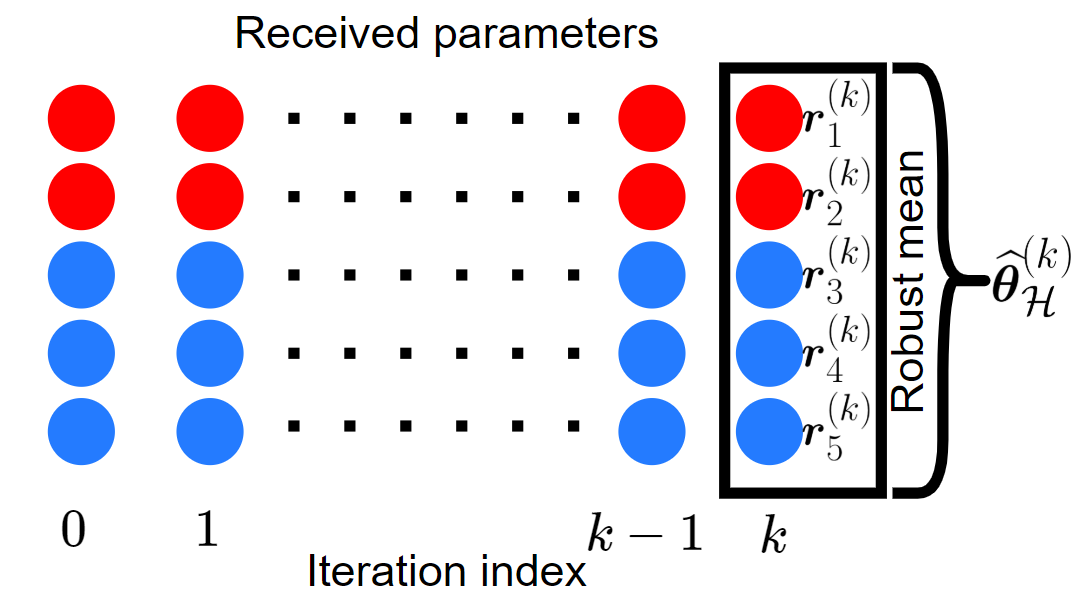}
    \caption{Robust mean estimation under static impersonation attack. Red/blue circles correspond to parameters received through compromised/trustworthy channels, respectively. In this example, there are $N=5$ agents and agents 1 and 2 are always communicating through compromised channels. At iteration $k$, the central coordinator computes the robust mean $\widehat{\prm}^{(k)}_{\cal H}$ of the received parameters $\{\bm r_i^{(k)}\}_{i\in[N]}$.}
    \label{fig:robust_pddra}
\end{figure}

We summarize the static impersonation attack resilient PD-DRA method in Algorithm~\ref{alg:robust_pddra}. The algorithm behaves similarly as Algorithm~\ref{alg:pddra} applied to \eqref{eq:pdp}, with the exception that the central coordinator is oblivious to ${\cal H}$, and it uses a robust mean estimator to find an approximate average for the signals sent through the trustworthy links, as illustrated in Figure~\ref{fig:robust_pddra}. This approximate value is used to compute the new price signals, and sent back to  agents.
In particular, the primal-dual updates are:
\begin{subequations}
\beq \label{eq:primal_pb}
\prm_i^{(k+1)} \hspace{-.2cm} = {\cal P}_{\Cset_i} \big( \prm_i^{(k)} \hspace{-.1cm} - \hspace{-.05cm} {\textstyle \frac{\gamma}{N}} \big( 
\widehat{\bm g}^{(k)}_{\cal H} + \grd_{\prm_i}\bt{f_i}( \prm_i^{(k)} ) + \upsilon \prm_i^{(k)} \big) \big),
\eeq
\beq \label{eq:dual_pb}
\lambda_t^{(k+1)} = \big[ \lambda_t^{(k)} + \gamma \big( \overline{g}_t ((1-\alpha_1) \widehat{\prm}_{\cal H}^{(k)} ) - \upsilon \lambda_t^{(k)} \big) \big]_+. 
\eeq
\end{subequations}
\btt{We note that the update rule in \eqref{eq:primal_pb} is valid for agents in set ${\cal H}$, because the gradients of the Lagrangian are defined only for those agents in \eqref{eq:grd_prm_robust}. The agents in set $\cal A$ may or may not use the same update rule, however, this does not have any impact on the algorithm as they can never communicate their true parameters to the central coordinator.}
\begin{lemma} \label{lem:perturbedgradients1}
Algorithm~\ref{alg:robust_pddra} is  a primal-dual algorithm \cite{koshal2011multiuser} for \eqref{eq:pdp} with perturbed gradients:\begin{subequations} \label{eq:perturbed}
\begin{align} 
\widehat{\bm g}_{ \prm }^{(k)} & = \grd_{\prm} \overline{\cal L}_\upsilon ( \prm^{(k)}; \bm{\lambda}^{(k)} ; {\cal H} ) + {\bm e}_{\prm}^{(k)} , \label{eq:per_prm} \\[.1cm]
\widehat{\bm g}_{ \bm{\lambda} }^{(k)} & =  \grd_{\bm{\lambda}} \overline{\cal L}_\upsilon (  \prm^{(k)} ; \bm{\lambda}^{(k)} ; {\cal H} ) + {\bm e}_{\bm{\lambda}}^{(k)},
\label{eq:per_lam}
\end{align}
\end{subequations}
where we have used concatenated variable as $\prm = (\{\prm_i\}_{i\in \cal{H}})$. Under Assumption \ref{ass:bdd} and assuming that $\lambda_t^{(k)} \leq \overline{\lambda}$ for all $k$, we have:
\beq\label{eq:bd_prm}
\begin{aligned}
\|  {\bm e}_{\prm}^{(k)} \| \leq& (1-\alpha_1) \overline{\lambda}L T  \| 
\widehat{\prm}_{\cal H}^{(k)} - \overline{\prm}_{\cal H}^{(k)} \|\\
&+\frac{|\mathcal H|-(1-\alpha_1)N}{|\mathcal H|}\overline{\lambda}B T,
\end{aligned}
\eeq
\beq \|{\bm e}_{\bm{\lambda}}^{(k)}\| \leq(1-\alpha_1) BT \| \widehat{\prm}_{\cal H}^{(k)} - \overline{\prm}_{\cal H}^{(k)} \|. \label{eq:bd_lambda}
\eeq
\end{lemma}
The proof can be found in Appendix~\ref{app:permanentperturbedgrads}. The assumption $\lambda_t^{(k)} \leq \overline{\lambda}$ can be guaranteed since $\overline{g}_t ( (1-\alpha_1) \widehat{\prm}_{\cal H}^{(k)} )$ is bounded\bt{, which is proven in Appendix~\ref{app:boundedgradient}.} 
Furthermore, the performance analysis for the median based estimator 
shows that 
\beq \label{eq:informal_bd}
\| \widehat{\prm}_{\cal H}^{(k)} - \overline{\prm}_{\cal H}^{(k)} \| = {\cal O}( \alpha_1 R \sqrt{d})
\eeq
when $\alpha_1$ is small.
Finally, based on Lemma~\ref{lem:perturbedgradients1}, we can analyze the convergence of Algorithm~\ref{alg:robust_pddra}. Let $\widehat{\bm z}^\star = (\widehat{\prm}^\star, \widehat{\bm{\lambda}}^\star)$
be a saddle point of \eqref{eq:pdp} and define 
\beq
\overline{\bm{\Phi}}( {\bm z}^{(k)} ) 
\eqdef \left( \begin{array}{c} 
\grd_{\prm} \overline{\cal L}_{\upsilon} ( \{ \prm_i^{(k)} \}_{i \in \cal{H}}, \bm{\lambda}^{(k)}; {\cal H} ) \\
- \grd_{\bm{\lambda}} \overline{\cal L}_{\upsilon} ( \{ \prm_i^{(k)} \}_{i \in \cal{H}}, \bm{\lambda}^{(k)} ; {\cal H} ) 
\end{array} \right).
\eeq 
We are ready to present our main result for static attacks.

\begin{theorem} \label{thm:main}
Assume the map $\overline{\bm{\Phi}}( {\bm z}^{(k)} )$ is $L_{\Phi}$-Lipschitz continuous. For Algorithm~\ref{alg:robust_pddra}, for all $k \geq 0$ it holds:
\beq \label{eq:rate_E} \begin{split}
 \| {\bm z}^{(k+1)} - \widehat{\bm z}^\star \|^2 \leq& \big( 1 - \gamma \upsilon' + 2 \gamma^2 L_{\Phi}^2 \big) \| {\bm z}^{(k)} - \widehat{\bm z}^\star \|^2 \\
&+ \left( \frac{4 \gamma}{\upsilon'} + 2 \gamma^2  \right) E_k.  \end{split}
\eeq
where $\upsilon'\eqdef(1-\alpha_1)\upsilon$ and $E_k \eqdef \| {\bm e}_{\prm}^{(k)} \|^2 + \| {\bm e}_{\bm{\lambda}}^{(k)} \|^2$ is the total perturbation at iteration $k$. Moreover, if we choose 
$\gamma < {\upsilon'} / {2L_{\Phi}^2}$
and $E_k$ is upper bounded by $\overline{E}$ for all $k$, then
\beq \label{eq:upperbd_E}
\limsup_{k \rightarrow \infty} \| {\bm z}^{(k)} - \widehat{\bm z}^\star \|^2 \leq 
\frac{ \frac{4}{\upsilon'} + 2 \gamma }{ \upsilon' - 2 \gamma L_{\Phi}^2} \overline{E}.
\eeq
\end{theorem}
The proof can be found in Appendix~\ref{app:permanentthm}. Combining with \eqref{eq:informal_bd} shows that the resilient PD-DRA method converges to a ${\cal O}(\alpha_1^2 R^2 d)$ neighborhood of the saddle point of \eqref{eq:pdp}.
Moreover, it shows that the convergence rate to the neighborhood is linear, which is similar to the classical analysis in \cite{koshal2011multiuser}.


\subsection{Dynamic Impersonation Attack}\label{subsec:infrequent}

Under this attack scenario, the set of agents communicating through compromised channels is dynamically changing with iterations. We make the following assumption on how frequently each agent's communications are compromised:
\begin{assumption}\label{ass:attackedfraction}
Let $m$ be a fixed window size and $\alpha_2<0.5$ be a known upper bound on how frequent an agent communicates through a compromised channel. Then, for all $k\geq m-1$ and for all agents $i\in[N]$, among the received parameters $\{\bm r_i^{(k-\ell)}\}_{\ell=0}^{m-1}$ at most $\alpha_2 m$ are sent through compromised channels.
\end{assumption}

It is important to recall that the dynamic attack scenario does not generalize
the static attack scenario and there is a significant distinction between the two. The static attack scenario assumes that a fixed set of agents' communications are \emph{permanently} compromised. \bt{It may occur when when the attacker compromises set of communication channels and those channels are assigned to the agents via a static channel allocation scheme.}

On the contrary, for the dynamic attack scenario, each user's communications are vulnerable to attacks for \emph{at most} a given $\alpha_2$ fraction of iterations over a window of size $m$ under Assumption~\ref{ass:attackedfraction}, and hence each agent is able to communicate reliably with the central coordinator at some iterations. \bt{This scheme may occur when the attacker compromises a fixed set of communication channels (same as the static scenario), however, the channels are assigned to the agents via a dynamic channel allocation scheme (e.g., do a round-robin channel allocation. If there are $m$ communication channels out of which $\alpha_2m$ are compromised, assigning channels dynamically in a cyclic way to the agents ensures that over a window of $m$, every agent has sent $\alpha_2m$ compromised messages). Although the attacker behaves the same way, we can simulate both scenarios by static/dynamic channel allocation.  In cyber-physical systems, such dynamic allocation schemes are commonly used (e.g., dynamic IP assignment to be protected from hackers).} 

Interestingly, it is possible to develop an algorithm that converges to the optimal solution of Problem \eqref{eq:pdreg}. 
The intuition behind is that the received parameters over a long period of time contain a fraction of trustworthy information that can extracted by the algorithm to perform faithful computations.

\begin{figure}[t]
    \centering
    \includegraphics[width=0.5\textwidth]{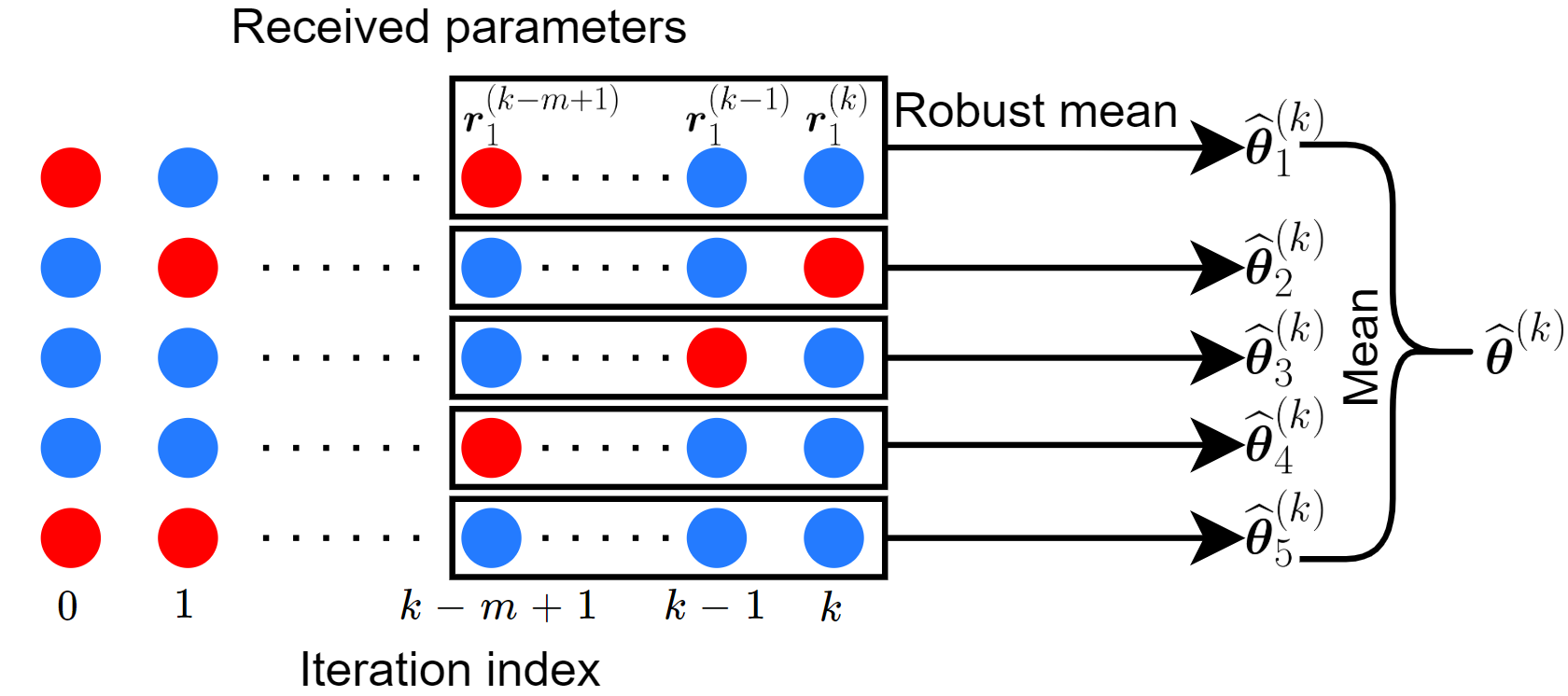}
    \caption{Robust mean estimation under dynamic impersonation attacks. Red/blue circles correspond to parameters received through compromised/trustworthy channels, respectively. In this example, there are $N=5$ agents and the set of agents communicating through compromised channels is changing at every iteration. At iteration $k$, the central coordinator computes the robust mean $\widehat{\prm}_i^{(k)}$ of the received parameters $\{\bm r_i^{(k-\ell)}\}_{\ell=0}^{m-1}$ for all agents $i\in[N]$. Then, computes the naive average of $\{\widehat{\prm}_i^{(k)}\}_{i=1}^N$ to get the average parameter $\widehat{\prm}^{(k)}$.}
    \label{fig:averaging_pddra}
\end{figure}

Our algorithm is similar in nature to an \emph{averaging gradient} scheme where the primal-dual updates utilize the averages of time delayed gradients. Furthermore, the scheme is combined with the \emph{robust mean} estimator developed in Sec.~\ref{sec:robustmean} to approximate the averages of outdated gradients, as illustrated in Figure~\ref{fig:averaging_pddra}. Specifically, the central coordinator chooses a window size of $m$. For any iteration $k\geq m-1$, instead of using $\bm r_i^{(k)}$ for computing the average parameter $\overline{\prm}^{(k)}$ and the gradients, the central coordinator computes the robust mean $\widehat{\prm}_i^{(k)}$ from the received parameters $\{\bm r_i^{(k-\ell)}\}_{\ell=0}^{m-1}$ using the median-based mean estimator \eqref{eq:median} for all agents $i\in[N]$, applied on the sequence of historical received parameters. Note that we have replaced $\alpha_1$ by $\alpha_2$, $N$ by $m$ in this application. It then uses $\widehat{\prm}^{(k)}\eqdef\frac{1}{N}\sum_{i=1}^N\widehat{\prm}_i^{(k)}$ for computation of the primal-dual updates.  

We summarize our robust averaging PD-DRA method in  Algorithm~\ref{alg:jamming_alg}. The primal-dual updates are described by:
\begin{subequations}\label{eq:updates_jam}
\beq \label{eq:prm_update_jam}
\prm_i^{(k+1)} \hspace{-.2cm} = {\cal P}_{\Cset_i} \big( \prm_i^{(k)} \hspace{-.1cm} - \hspace{-.05cm} {\textstyle \frac{\gamma}{N}} \big( 
\widehat{\bm g}^{(k)} + \grd_{\prm_i} \bt{f_i}( \prm_i^{(k)} ) + \upsilon \prm_i^{(k)} \big) \big),
\eeq
\beq \label{eq:dl_update_jam}
\lambda_t^{(k+1)} = \big[ \lambda_t^{(k)} + \gamma \big( g_t ( {\textstyle \widehat{\prm}^{(k)} ) - \upsilon \lambda_t^{(k)} \big) \big]_+}. 
\eeq
\end{subequations}
\begin{lemma}\label{lemma:perturbedgradients2}
Algorithm \ref{alg:jamming_alg} is a primal-dual algorithm for \eqref{eq:pdreg} with perturbed gradients:
\begin{subequations}\label{eq:perturbgradients2}
\begin{align}
    \label{eq:thetaperturb2}&\boldsymbol{\widehat{g}}^{(k)}_{\prm}=\nabla_{\prm}\mathcal{L}_\upsilon(\prm^{(k)};\boldsymbol{\lambda}^{(k)})+\boldsymbol{e}_{\prm}^{(k)},\\
    \label{eq:lambdaperturb2}&\boldsymbol{\widehat{g}}^{(k)}_{\boldsymbol{\lambda}}=\nabla_{\boldsymbol{\lambda}} \mathcal{L}_\upsilon(\prm^{(k)};\boldsymbol{\lambda}^{(k)})+\boldsymbol{e}_{\boldsymbol{\lambda}}^{(k)},
    \end{align}
\end{subequations}
where we have used concatenated variable as $\prm=(\{\prm_i\}_{i\in N})$. Under Assumption \ref{ass:bdd} and assuming that $\lambda_t^{(k)}\leq\overline{\lambda}$ for all k, we have:
\begin{subequations}\label{eq:errorgrad2}
\begin{align}
    \label{eq:errorthetagrad2}&\|\boldsymbol{e}^{(k)}_{\prm}\|\leq \frac{\overline{\lambda}LT}{N}\sum_{i=1}^N\|\prm_i^{(k)}-\widehat{\prm}{}_i^{(k)}\|,\\
    \label{eq:errorlambdagrad2}&\|\boldsymbol{e}^{(k)}_{\boldsymbol{\lambda}}\|\leq \frac{BT}{N}\sum_{i=1}^N\|\prm_i^{(k)}-\widehat{\prm}{}_i^{(k)}\|.
    \end{align}
\end{subequations}
\end{lemma}

\begin{algorithm}[t]
	\caption{Averaging PD-DRA Algorithm}\label{alg:jamming_alg}
	\begin{algorithmic}[1]
		\STATE \textbf{Input}: Each agent has initial state $ \prm_i^{(0)}$.
		\FOR {$k=0,1,\dots,m-2$}
		\STATE Apply basic PD-DRA (Run Algorithm~\ref{alg:pddra}).
		\ENDFOR
		\FOR {$k=m-1,m,\dots$}
		\STATE \label{av_step:accum0}\emph{(At the Central Coordinator)}:
		\begin{enumerate}[label=(\alph*)]
			\item Receives $\{{\bm r}_i^{(k)}\}_{i=1}^N$, see~\eqref{eq:signals}, from agents.
			\item \label{av_step:averaging}For all agents $i=1,\dots, N$, computes robust mean  $\widehat{\prm}_i^{(k)}$ of $\{\bm r_i^{(k-\ell)}\}_{\ell=0}^{m-1}$ using the estimator \eqref{eq:median} with parameters $\alpha_1\rightarrow \alpha_2$, $N\rightarrow m$.
			\item Computes $\widehat{\prm}^{(k)}\eqdef\frac{1}{N}\sum_{i=1}^N\widehat{\prm}^{(k)}_i$.
			\item Broadcasts the vector $\widehat{\bm g}^{(k)} \eqdef \sum_{t=1}^T \lambda_t^{(k)} \nabla_{\prm} g_t (\widehat{\prm}^{(k)}) $ to  agents.
			\item Computes the update for $\bm{\lambda}^{(k+1)}$  with \eqref{eq:dl_update_jam}. 
		\end{enumerate}
		\STATE \emph{(At each agent $i$)}:
		\begin{enumerate}[label=(\alph*)]
			\item Agent receives $\widehat{\bm g}^{(k)}$.
			\item Agent computes update for $\prm_i^{(k+1)}$ with \eqref{eq:prm_update_jam}.

		\end{enumerate}
		\ENDFOR
	\end{algorithmic} 
\end{algorithm}

\noindent The proof can be found in Appendix~\ref{app:perturbedgradientsinfreq}. The assumption $\lambda_t^{(k)} \leq \overline{\lambda}$ can be guaranteed since $g_t(\widehat{\prm}^{(k)})$ is bounded\bt{, which is proven in Appendix~\ref{app:boundedgradient}}. Let ${\bm z}^{(k)} \eqdef ( \{ \prm_i^{(k)} \}_{i=1}^N, \bm{\lambda}^{(k)} )$
be the primal-dual variable at the  $k$th iteration and define the mapping $\bm{\Phi}( {\bm z}^{(k)} )$ as in \eqref{eq:mapping}. 
We observe that the algorithm's behavior is similar to the incremental aggregated gradient method in \cite{iag_blatt,iag_metasupablo,iag_Tseng2014}.
The following Lemma, which is inspired by \cite{iag_blatt,iag_metasupablo,iag_Tseng2014}, upper bounds the perturbation in the gradients in \eqref{eq:errorgrad2} by the maximum optimality gap in a finite window of size $2m-1$:
\begin{lemma}\label{lemma:perturbationupperbound}
Assume the map $\boldsymbol{\Phi}(\prmdl^{(k)})$ is $L_{\Phi}$-Lipschitz continuous. Let $E_k:=\|\boldsymbol{e}^{(k)}_{\prm}\|^2+\|\boldsymbol{e}^{(k)}_{\boldsymbol{\lambda}}\|^2$. Then, for all $k\geq2(m-1)$ we have:
\begin{equation}
E_k\leq \gamma^2 \overline{C}\underset{0\leq \ell\leq2(m-1)}{\max}\|\prmdl^{(k-\ell)}-\prmdl^\star\|^2,
\end{equation}
where 
\beq \label{eq:constantcbar}
\begin{aligned}
   \overline{C}=&\left(\frac{T^2(\overline{\lambda}^2L^2+B^2)}{N}\right)\times\left(\frac{1}{L_\Phi}+(1+\sqrt{d})\overline{\lambda}LT\right)^2\\&\times \left(\frac{1+C_\alpha}{1-\alpha_2}+C_\alpha\right)^2\times(m-1)^2,  
\end{aligned}
\eeq
and
\begin{equation*}
    C_\alpha=\frac{2\alpha_2}{1-\alpha_2}\left(1+\sqrt{\frac{(1-\alpha_2)^2}{1-2\alpha_2}}\right)\sqrt{d}.
\end{equation*}
\end{lemma}
\noindent The proof can be found in Appendix~\ref{app:perturbationbound}. Using on Lemmas~\ref{lemma:perturbedgradients2} and \ref{lemma:perturbationupperbound}, we can analyze the converge of Algorithm~\ref{alg:jamming_alg}:
\begin{theorem}\label{thm:convergence}
Assume the map $\boldsymbol{\Phi}(\prmdl^{(k)})$ is $L_{\Phi}$-Lipschitz continuous. For Algorithm~\ref{alg:jamming_alg}, for all $k\geq 2(m-1)$ it holds that:
\begin{equation}\label{eq:thm1eq1}
\begin{aligned}
    &\|\prmdl^{(k+1)}-\prmdl^\star\|\btt{^2}\leq (1-\gamma \upsilon+2\gamma^2L_\Phi^2)\|\prmdl^{(k)}-\prmdl^\star\|^2\\
    &+\left(\frac{4\gamma}{\upsilon}+2\gamma^2\right)\gamma^2\overline{C}\underset{0\leq \ell\leq2(m-1)}{\max}\|\prmdl^{(k-\ell)}-\prmdl^\star\|^2.
\end{aligned}
\end{equation}
Moreover, if we choose $\gamma$ sufficiently small such that it satisfies
$$\upsilon-2\gamma L_\Phi^2-\frac{4\overline{C}\gamma^2}{\upsilon}-2\overline{C}\gamma^3>0,$$
then:
\begin{equation}\label{eq:geometricconvergence}
    \|\prmdl^{(k)}-\prmdl^\star\|^2\leq\rho^{k-2(m-1)}\|\prmdl^{(2(m-1))}-\prmdl^\star\|^2,
\end{equation}
and
\begin{equation}
    \label{eq:thm1eq2}
    \underset{k\rightarrow\infty}{\lim}\|\prmdl^{(k)}-\prmdl^\star\|^2=0,
\end{equation}
where $\rho=(1-\gamma \upsilon+2\gamma^2L_\Phi^2+\frac{4\overline{C}\gamma^3}{\upsilon}+2\overline{C}\gamma^4)^\frac{1}{1+2(m-1)}$. 
\end{theorem}
The proof can be found in Appendix~\ref{app:infrequentthm}. Theorem~\ref{thm:convergence} shows that the robust averaging PD-DRA method converges \emph{geometrically} to the optimal solution of \eqref{eq:pdreg} under said assumptions. 

\subsection{Remarks}

A few remarks highlighting design criteria to be explored in practical implementations are in order:
\begin{itemize}[leftmargin=5mm]
    \item Theorem~\ref{thm:main} illustrates a trade-off in the choice of the step size $\gamma$ between convergence speed and accuracy. In particular, \eqref{eq:rate_E} shows that the rate of convergence  factor $1 - \gamma \upsilon + 2 \gamma^2 L_\Phi^2$ can be minimized by setting $\gamma = \upsilon / (4L_\Phi^2)$. 
    Meanwhile, the asymptotic upper bound in \eqref{eq:upperbd_E} is increasing with $\gamma$ and it can be minimized by setting $\gamma \rightarrow 0$.
    \item Theorem~\ref{thm:convergence} illustrates a trade-off between 
    the window size $m$ and the convergence rate. Observe that increasing the window size $m$ decreases the rate of convergence by increasing $\rho$ (Equations \eqref{eq:constantcbar} and \eqref{eq:geometricconvergence}). On the other hand, the likelihood that Assumption~\ref{ass:attackedfraction} holds true in a stochastic setting (e.g., channels being compromised with some probability) increases with a larger window size $m$.
    \item Under the dynamic impersonation attack scenario, the choice of $\alpha_2$ does not affect convergence accuracy to the saddle point of \eqref{eq:pdreg}, but only changes the convergence rate. As such, choosing the largest $\alpha_2$ such that $\alpha_2m=\lfloor\frac{m-1}{2}\rfloor$ (i.e., assuming maximum possible number of iterates received through compromised channels) makes the algorithm robustly applicable to all dynamic impersonation attack scenarios regardless of the frequency of the attack.
    \item In case the central coordinator can not identify the attack scenario as static or dynamic impersonation (or the attack can be a mixture of both), a mixture of both Algorithms~\ref{alg:robust_pddra} and \ref{alg:jamming_alg} can be applied. In particular, this can be done by adding Step~\ref{av_step:accum0}\ref{av_step:averaging} of Algorithm~\ref{alg:jamming_alg} before Step~\ref{r_step:accum0}\ref{r_step:robustmean} of Algorithm~\ref{alg:robust_pddra}, and applying the rest of the Algorithm~\ref{alg:robust_pddra} as it is. The central coordinator first computes robust parameters $\widehat{\prm}_i^{(k)}$ by computing the robust mean of $\{\bm r_i^{(k-\ell)}\}_{\ell=0}^{m-1}$ for all agents, and then computes the robust mean of  $\{\widehat{\prm}_i^{(k)}\}_{i=1}^N$. This effectively makes Algorithm~\ref{alg:robust_pddra} robust to possible dynamic impersonation attacks on uplink channels that are thought to be trustworthy for all iterations.
\end{itemize}

\section{Numerical Study}\label{sec:numerical}

In this section, we demonstrate the performance of our methods and verify our theoretical claims by applying our algorithms for: 1) an electric vehicle (EV) charging coordinator under static impersonation attack, 2) an electric vehicle charging coordinator under dynamic impersonation attack, and  3) a power distribution network with flexible demand under dynamic impersonation attack. The EV charging coordinator problem resembles classic network utility maximization problems such as those studied in communication networks whereas the power distribution network problem has more nuisances that we will discuss next. \bt{To solve the convex optimization problems in order to get the optimal solutions, we used CVX, a package for specifying and solving convex programs\cite{cvx}.}

\subsection{Electric Vehicle (EV) Charging Facility}

In this study, the aim is to optimize EV charging demand over time. We consider multiple EVs receiving charge under the same local feeder/transformer. Each agent (or EV owner) has different utility of charging at different times. Hence, at a given time period, it is desired to charge those EVs who have a higher utility (or less cost) for  that time period. This problem falls into the broad category of network utility maximization problems, which can be formulated as:
\begin{subequations}
\label{eq:evchargeoptimization}
\begin{align}
&\underset{ \prm_i \in \mathbb R_+^d,  \forall i}{\min} & &\bt{f}(\prm)=\frac{1}{N} \sum_{i=1}^N \bt{f_i} ( \prm_i )\\
& \text{subject to}
& & \frac{1}{N}\sum_{i=1}^N\prm_i\preceq \overline{\boldsymbol{e}},\\
& & & \label{eq:evconst2}\prm_i^{\min}\preceq \prm_i \preceq \prm_i^{\max},\quad\forall i,\\
& & & \label{eq:evconst3}\Theta_i^{\min}\leq \boldsymbol{1}^T\prm_i\leq\Theta_i^{\max},\quad\forall i,
\end{align}
\end{subequations}
where $N\times\overline{\boldsymbol{e}}\in \mathbb R^d$ is the vector of maximum available transformer capacity in all time periods \btt{and $\preceq$ denotes component-wise inequality between the vectors}. The available capacity changes with time of day as exogenous load on the transformer varies with time as well. The elements $\{\prm_{i,j}\}_{j=1}^d$ of the vector $\prm_i$ correspond to the electricity demand of the EV $i$ at time slots $j=1\dots d$. The constraint \eqref{eq:evconst2} restricts the amount an EV can charge at each time slot, whereas the constraint \eqref{eq:evconst3} bounds the total amount an EV can charge. For this study, we set the cost function to be:
\begin{equation}
    \label{eq:utilitycharging}
    \bt{f_i}(\prm)=-\sum_{j=1}^d\beta_{i,j}\log\prm_{i,j},
\end{equation}
where $\beta_{i,j}$ are generated randomly from a uniform distribution in $[0,1]$. 
We study this problem under both attack scenarios for $N=100$ EVs.

\begin{figure*}[t]
     \centering
     \hspace{-0.5cm}
     \begin{subfigure}[b]{0.25\textwidth}
         \centering
         \includegraphics[width=\textwidth]{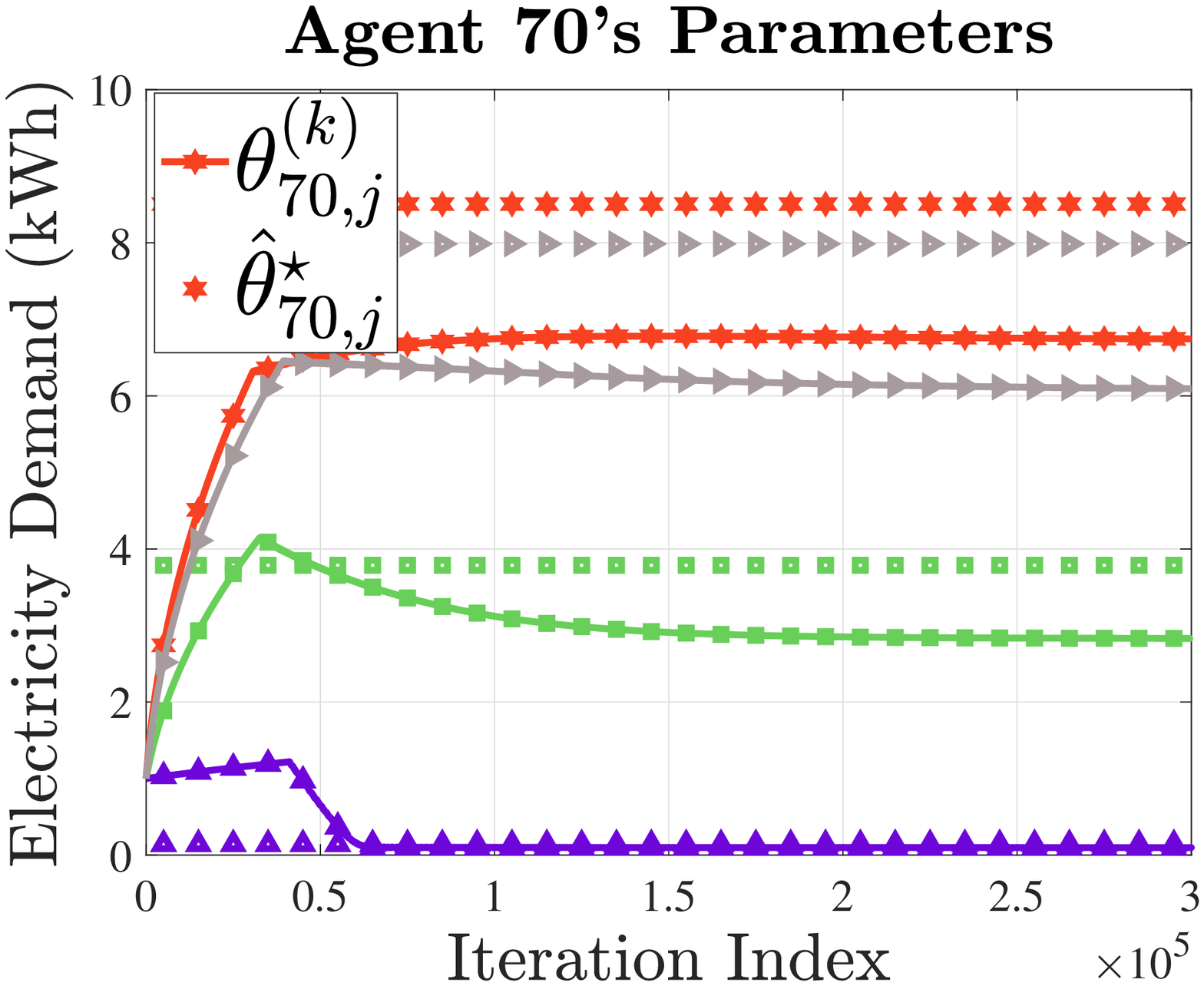}
         \caption{}
         \label{fig:evconvergence}
     \end{subfigure}
     \hspace{-0.5cm}
     \begin{subfigure}[b]{0.25\textwidth}
         \centering
         \includegraphics[width=\textwidth]{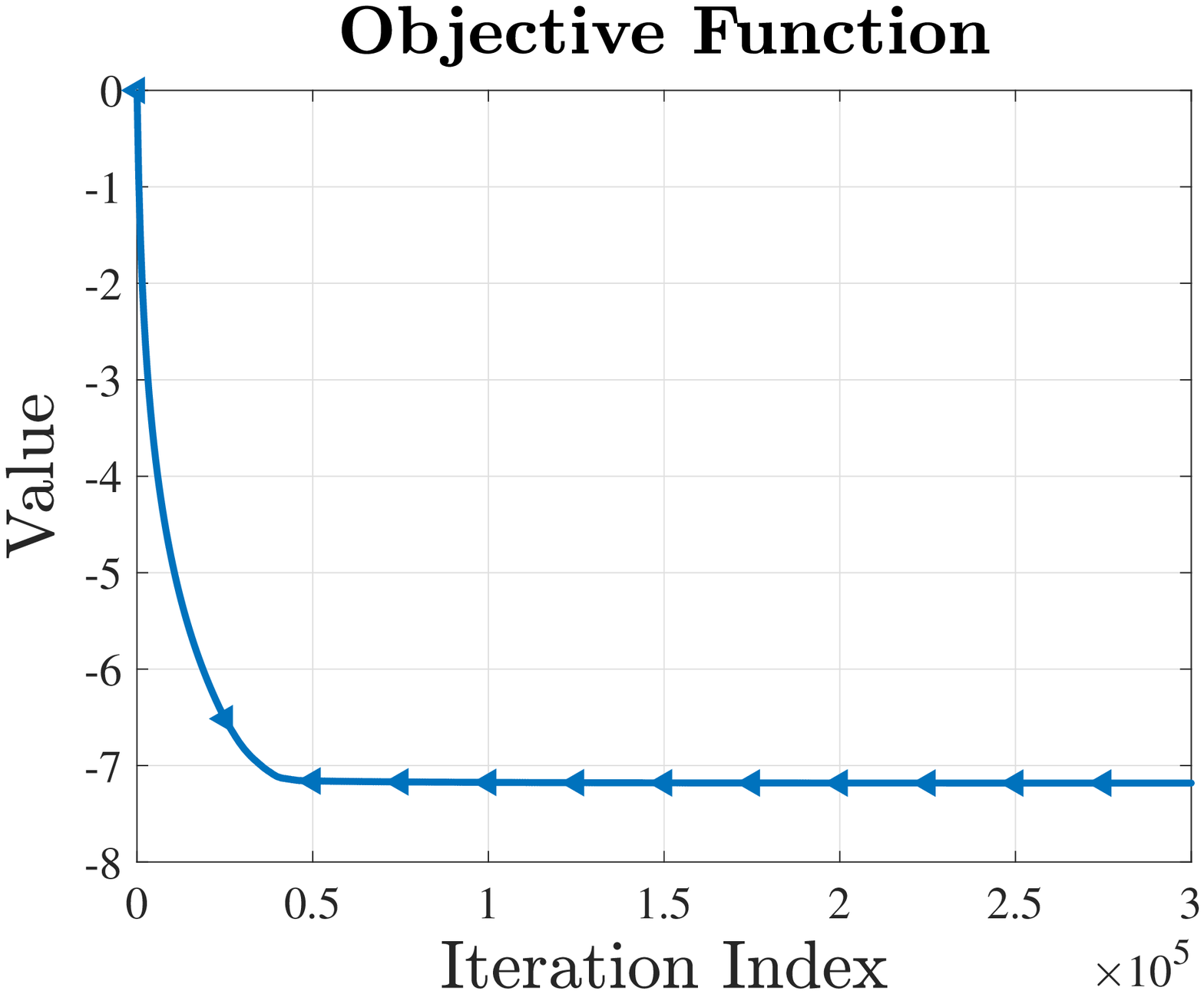}
         \caption{}
         \label{fig:evobjconst}
     \end{subfigure}
     \hspace{-0.5cm}
     \begin{subfigure}[b]{0.25\textwidth}
         \centering
         \includegraphics[width=\textwidth]{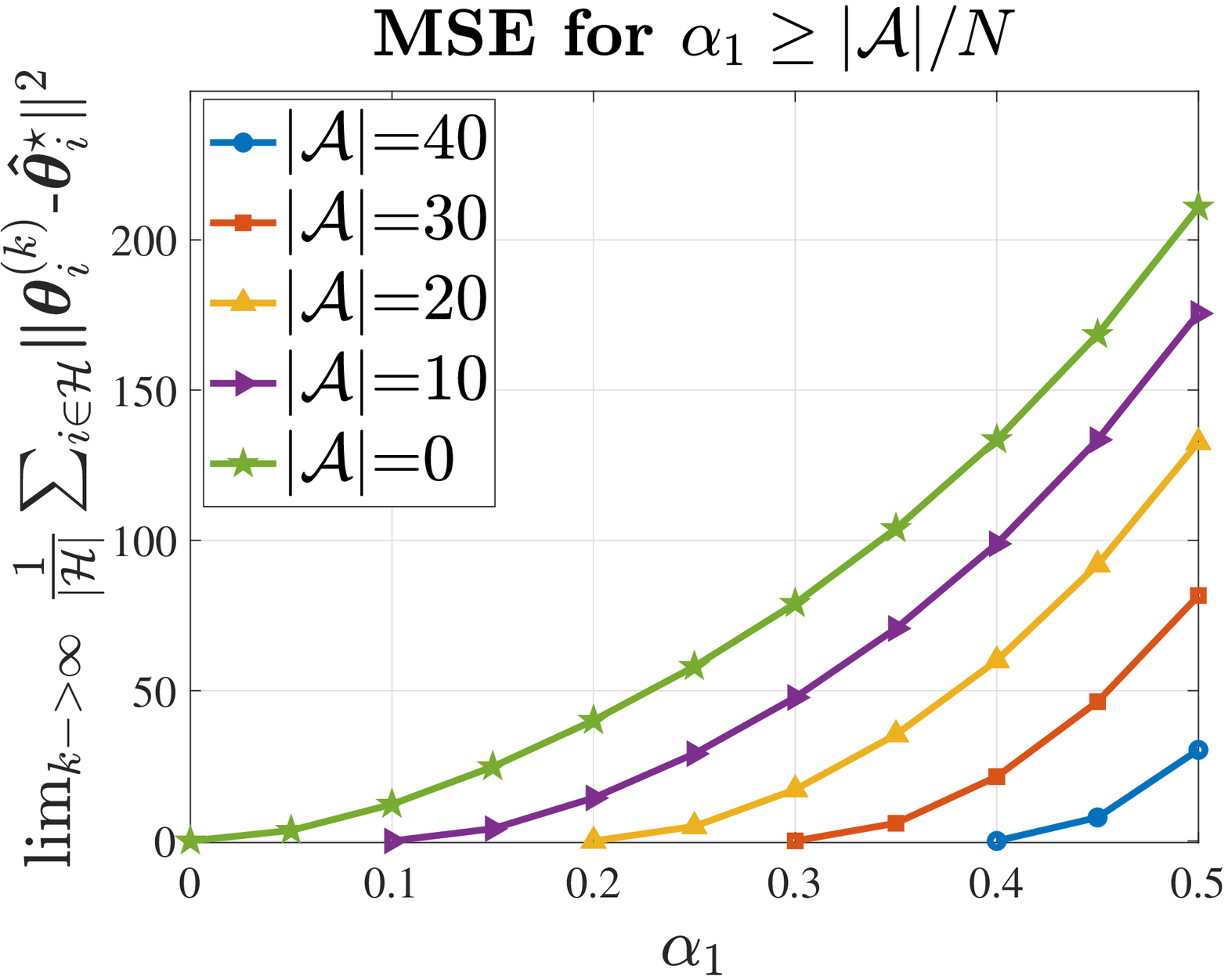}
         \caption{}
         \label{fig:evAvsalpha}
     \end{subfigure}
     \hspace{-0.5cm}
     \begin{subfigure}[b]{0.25\textwidth}
         \centering
         \includegraphics[width=\textwidth]{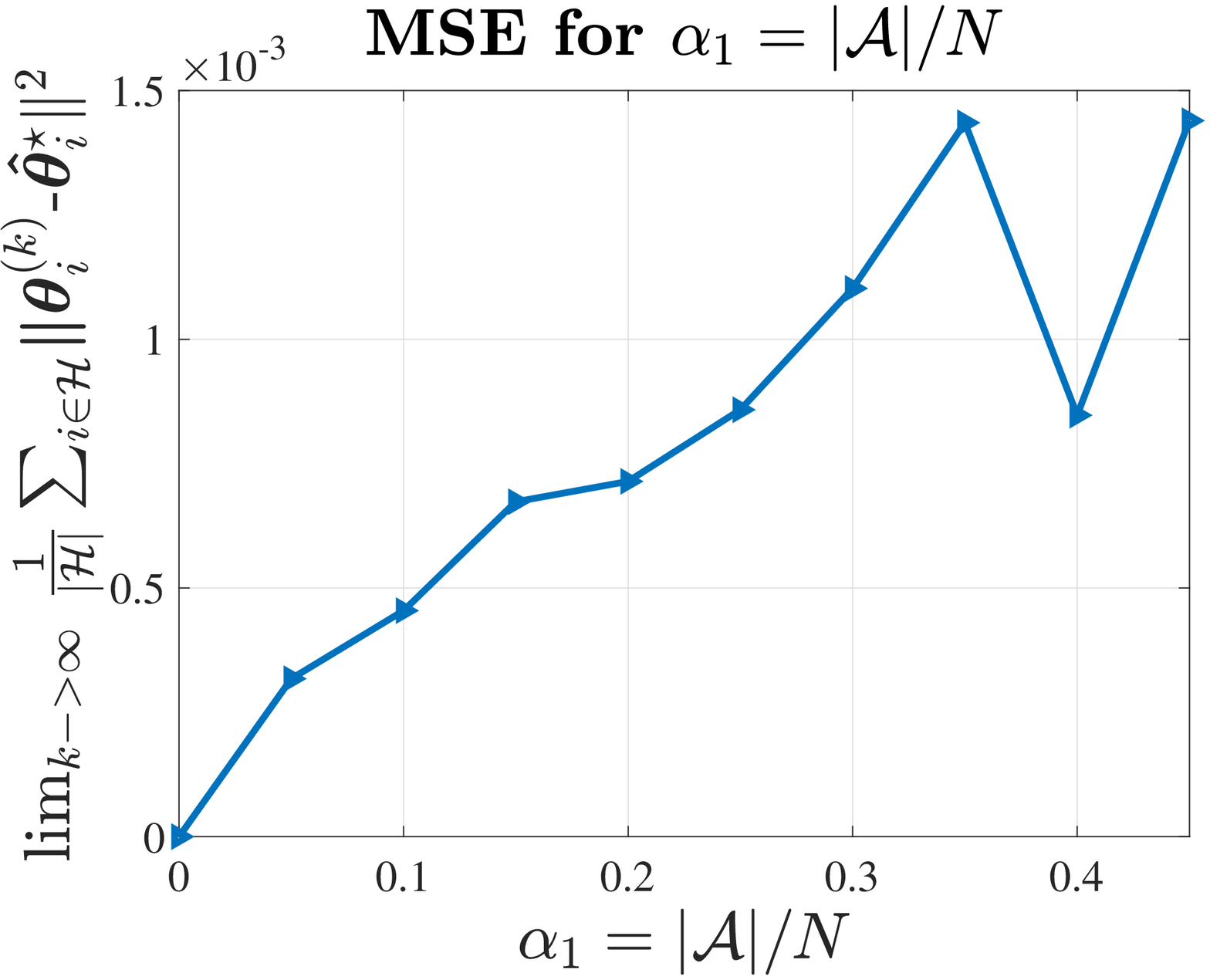}
         \caption{}
         \label{fig:evmsewithcorrectalpha}
     \end{subfigure}
        \caption{Numerical study results for optimal electric vehicle charging under static impersonation attack. (a) Optimal parameter of agent 70 converges to a neighborhood of the optimal solution of the robust optimization problem for $|{\cal A}|/N=0.2$ and $\alpha_1=0.3$, (b) The algorithm provides convergence of the objective function value, (c) Mean squared error for different number of compromised channels and different choices of upper bound $\alpha_1$, (d) Mean squared error when $\alpha_1=|{\cal A}|/N$.}
        \label{fig:evchargingresults}
              \vspace{-.4cm}
\end{figure*}

\subsubsection{Static impersonation attack} We simulated various static impersonation attack scenarios and ran Algorithm~\ref{alg:robust_pddra}. The results are displayed in Figure~\ref{fig:evchargingresults}.

In Figure~\ref{fig:evconvergence}, we plot agent 70's electricity demand for \btt{some} time periods, with $|{\cal A}|/N=0.2$ and $\alpha_1=0.3$. \bt{Each different color corresponds to a different dimension of the parameter vector (i.e., electricity demand for different time periods). A colored solid line corresponds to a dimension of the parameter vector iterates generated by the algorithm. A dashed line with the same \btt{marker and color} as a solid line is  the optimal value corresponding to that dimension of the parameter vector, which is the solution of the regularized robust optimization problem (formulated as \eqref{eq:pdp}) of \eqref{eq:evchargeoptimization}.} Observe that Algorithm~\ref{alg:robust_pddra} successfully provides convergence to a close neighborhood of the optimal solution of the regularized robust optimization problem. Furthermore, in Figure~\ref{fig:evobjconst} we show that the objective function value converges, as opposed to a non-resilient PD-DRA method that is shown to oscillate and violate the constraint in Figure~\ref{fig:nonresilientfailure}. Our robust optimization model on the other hand ensures there is no constraint violation.

\begin{figure}[t]
     \centering
     \hspace{-0.4cm}
     \begin{subfigure}[b]{0.25\textwidth}
         \centering
         \includegraphics[width=\textwidth]{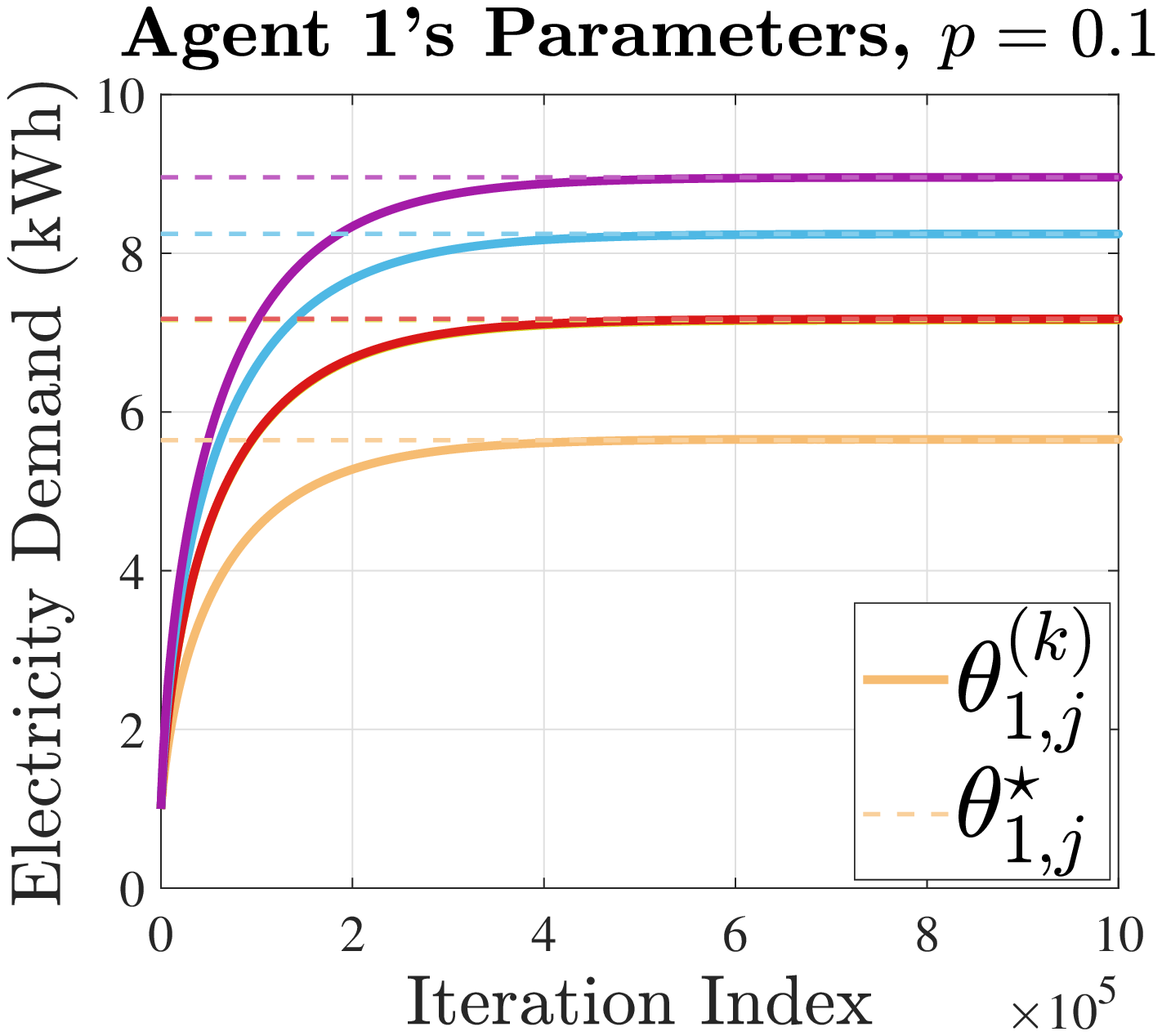}
         \caption{$p=0.1$, $m=20$, $\alpha_2=0.45$}
         \label{fig:evdynamicp01}
     \end{subfigure}
     \hspace{-0.2cm}
     \begin{subfigure}[b]{0.25\textwidth}
         \centering
         \includegraphics[width=\textwidth]{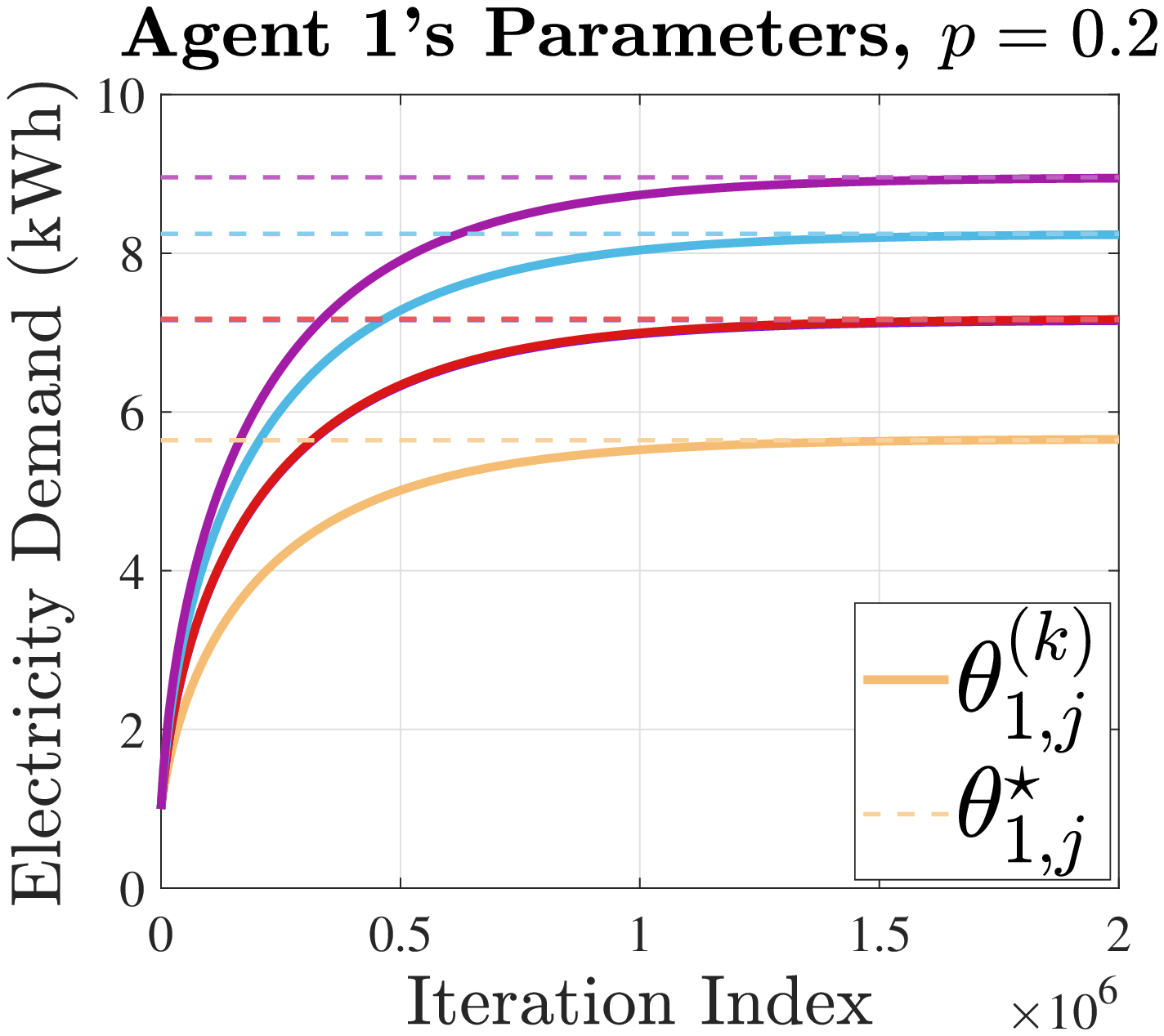}
         \caption{$p=0.2$, $m=100$, $\alpha_2=0.49$}
         \label{fig:evdynamicp02}
     \end{subfigure}
        \caption{Numerical study results demonstrating convergence of Algorithm~\ref{alg:jamming_alg} for optimal electric vehicle charging under two dynamic impersonation attack scenarios: (a) $p=0.1$, (b) $p=0.2$. Observe that the number of iterations it takes to converge for (b) is much larger than for (a).}
        \label{fig:evchargingdynamicresults}

\end{figure}

In Figure~\ref{fig:evAvsalpha}, we plot the mean squared error (MSE) in primal variables $\prm_i$ for different number of compromised channels $|\mathcal A|$ and different choices of $\alpha_1$, which is the upper bound on fraction of compromised links known by the central coordinator. The MSE is calculated by:
\begin{equation}
    \textnormal{MSE}=\underset{k\rightarrow\infty}{\lim}\frac{1}{|\mathcal{H}|}\sum_{i\in\mathcal{H}}\|\prm_i^{(k)}-\hat{\prm}^\star_i\|^2,
\end{equation}
where $\hat{\prm}^\star_i$ is the solution to $\eqref{eq:pdp}$ with $\alpha_1=|\mathcal{A}|/N$, \ie the solution to the regularized and robustified problem with the knowledge of the compromised channels. Naturally, the looser the upper bound, the larger the error, since it increases the amount of conservatism. Hence, having an accurate upper bound on fraction of compromised channels significantly improves the performance.

Finally, in Figure~\ref{fig:evmsewithcorrectalpha} we exhibit the efficacy of our approach with median-based mean estimation. We plot the mean squared error in primal variables, when the upper bound on $\alpha_1$ is tight, \ie $\alpha_1=|{\cal A}|/N$. The error tends to increase with $|{\cal A}|/N$, however, considering the magnitude, the error is negligible and we can conclude that the median-based mean estimator performs well.

\subsubsection{Dynamic impersonation attack} We simulated a dynamic impersonation attack scenario and ran Algorithm~\ref{alg:jamming_alg}. To simulate a dynamic impersonation attack, we assigned a probability $p$ for an uplink to be compromised at each iteration\footnote{Although a probabilistic scenario does not guarantee that Assumption~\ref{ass:attackedfraction} holds, with sufficiently large window size $m$ and $\alpha_2$, it holds with high probability at each iteration. Even though we do not study this scenario theoretically, our algorithm still performs well.}. For $p=0.1$, we picked a window size $m=20$ and $\alpha_2=0.45$, whereas for $p=0.2$, we picked a window size $m=100$ and $\alpha_2=0.49$. The results are displayed in Figure~\ref{fig:evchargingdynamicresults}. \bt{Each different color corresponds to a different dimension of the parameter vector. A colored solid line corresponds to a dimension of the parameter vector iterates generated by the algorithm. A dashed line with the same color as a solid line is  the optimal value corresponding to that dimension of the parameter vector, which is the solution of the regularized optimization problem (formulated as \eqref{eq:pdreg}) of \eqref{eq:evchargeoptimization}.}

In both scenarios, Algorithm~\ref{alg:jamming_alg} successfully provides convergence to the optimal solution of the regularized problem. Observe that for $p=0.2$, we chose a larger window size and a larger $\alpha_2$ in order to meet Assumption~\ref{ass:attackedfraction}. However, this restricts us to choose a smaller step size $\gamma$ as dictated by Theorem~\ref{thm:convergence} and in turn slower convergence. This highlights an important trade-off between robustness and convergence rate, where a larger window size $m$ and larger $\alpha_2$ makes the algorithm more robust while decreasing the convergence rate.
\subsection{Power Distribution Network}
\begin{figure}[t]
    \centering
    \includegraphics[width=.35\textwidth]{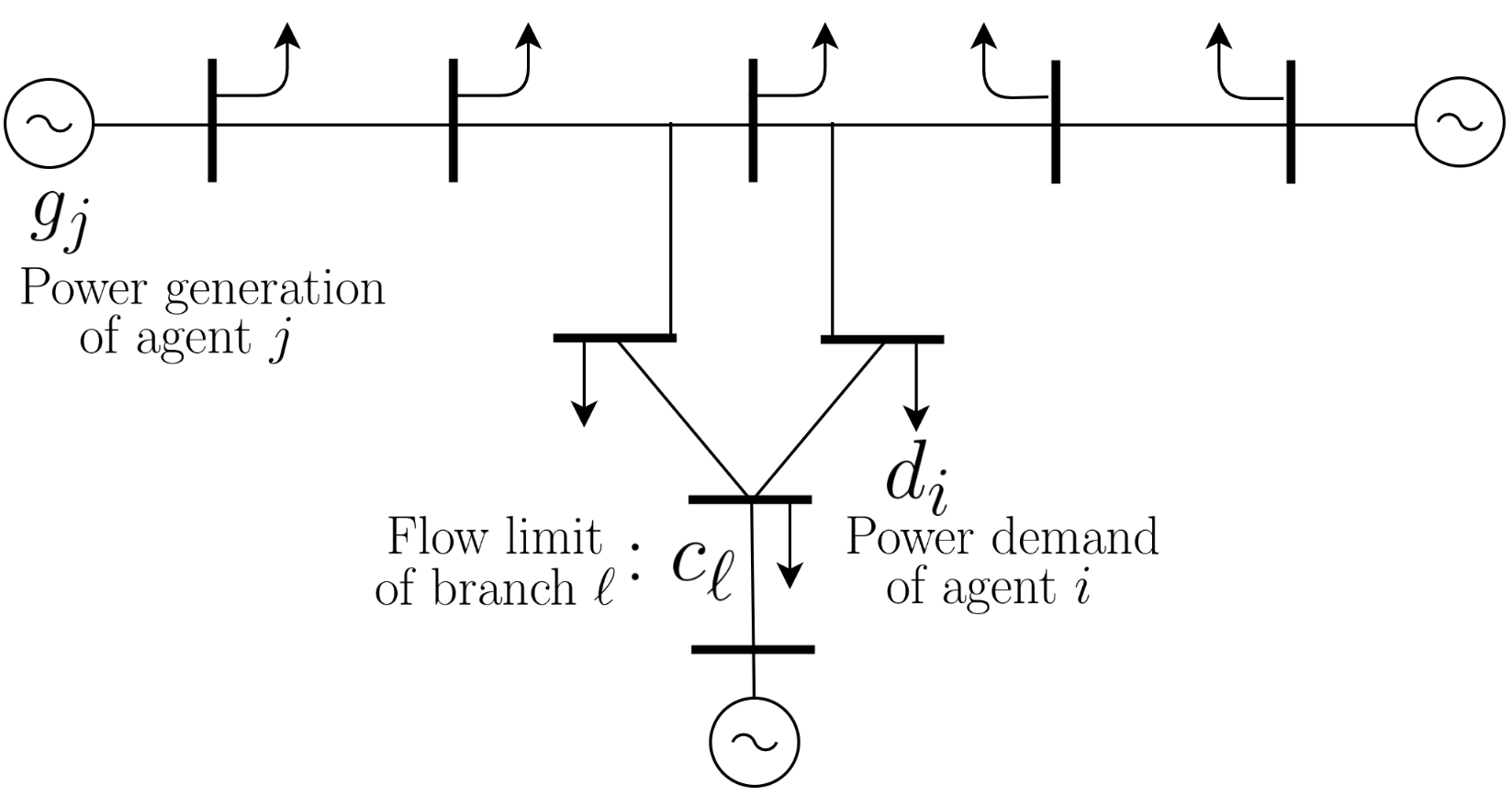}
    \caption{IEEE 9 bus system with 3 generators (supplies) represented by sources and 8 loads (demands) represented by arrows.}
    \label{fig:ieee9bus}
\end{figure}
We consider the IEEE $N=9$ bus system with $N_g=3$ generators and $N_{\ell}=8$ loads as shown in Figure~\ref{fig:ieee9bus}. The power network cost minimization problem can be stated as:
\begin{subequations}
\label{eq:poweroptimization}
\begin{align}
&\underset{ d_i, g_i \in \mathbb R^+}{\min} & &\bt{f}(\boldsymbol{d},\boldsymbol{g})= -\sum_{i=1}^{N_{\ell}}U_i ( d_i )+\sum_{i=1}^{N_g}C_i(g_i)\\
& \text{subject to}
\label{const:supplydemand}& & \boldsymbol{1}^T(\boldsymbol{d}-\boldsymbol{g})=0,\\
\label{const:flow}&&&\textnormal{\textbf{H}}(\boldsymbol{d}-\boldsymbol{g})\leq \boldsymbol{c},
\end{align}
\end{subequations}
where $\boldsymbol{d}=[d_1\dots d_N]^T$ and $\boldsymbol{g}=[g_1\dots g_N]^T$ are the vectors of load and generation at each node, respectively ($d_i=0$ for nodes without load and $g_j=0$ for nodes without generators). The first constraint \eqref{const:supplydemand} ensures the power supply is equal to the demand, and the second constraint \eqref{const:flow} is the power flow constraint limiting the power flow on each branch.

\begin{figure}[t]
     \centering
     \begin{subfigure}[b]{0.5\textwidth}
         \centering
         \includegraphics[width=\textwidth]{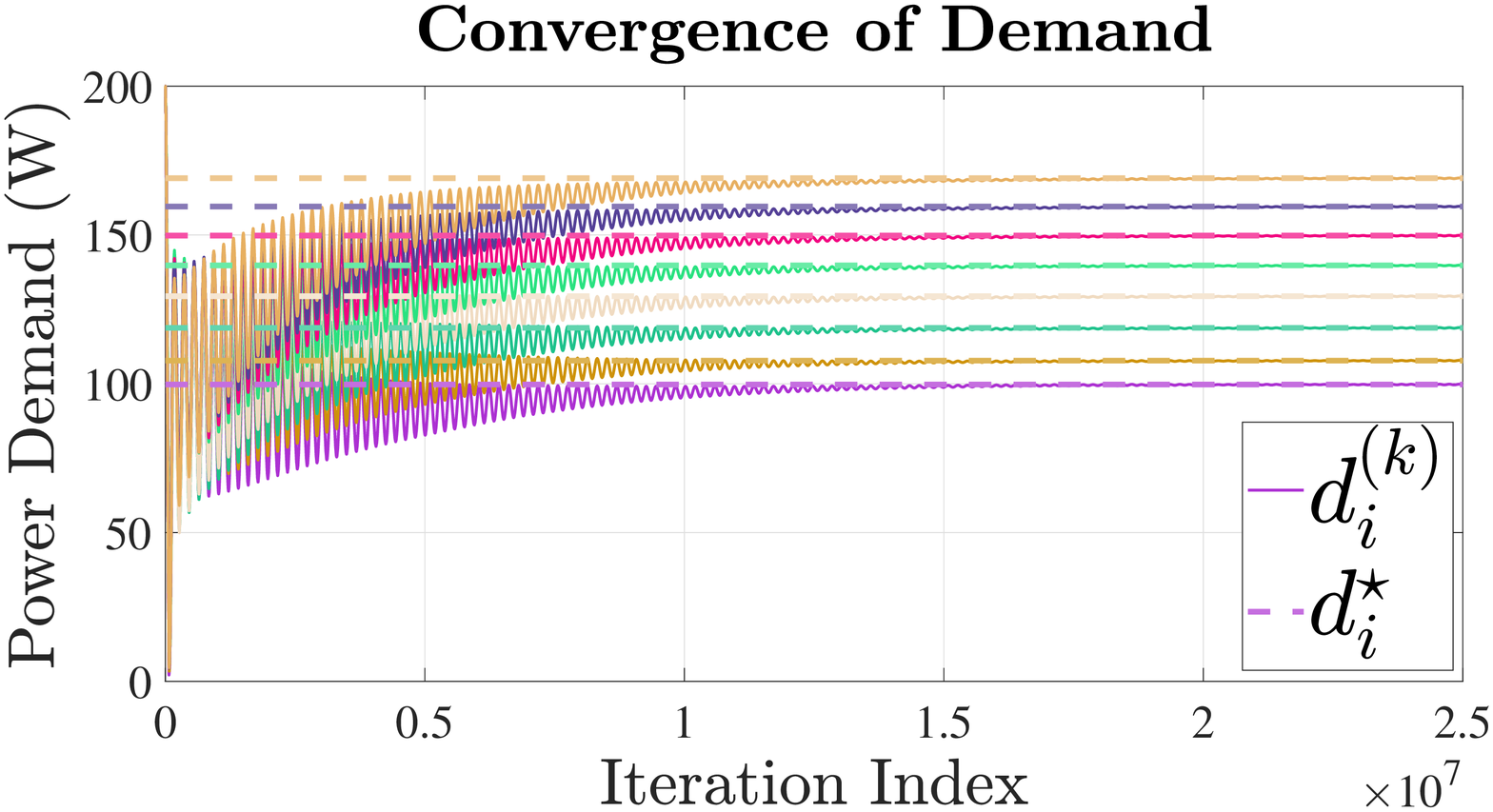}
         \caption{}
         \label{fig:powernetworkdemand}
     \end{subfigure}
     \begin{subfigure}[b]{0.5\textwidth}
         \centering
         \includegraphics[width=\textwidth]{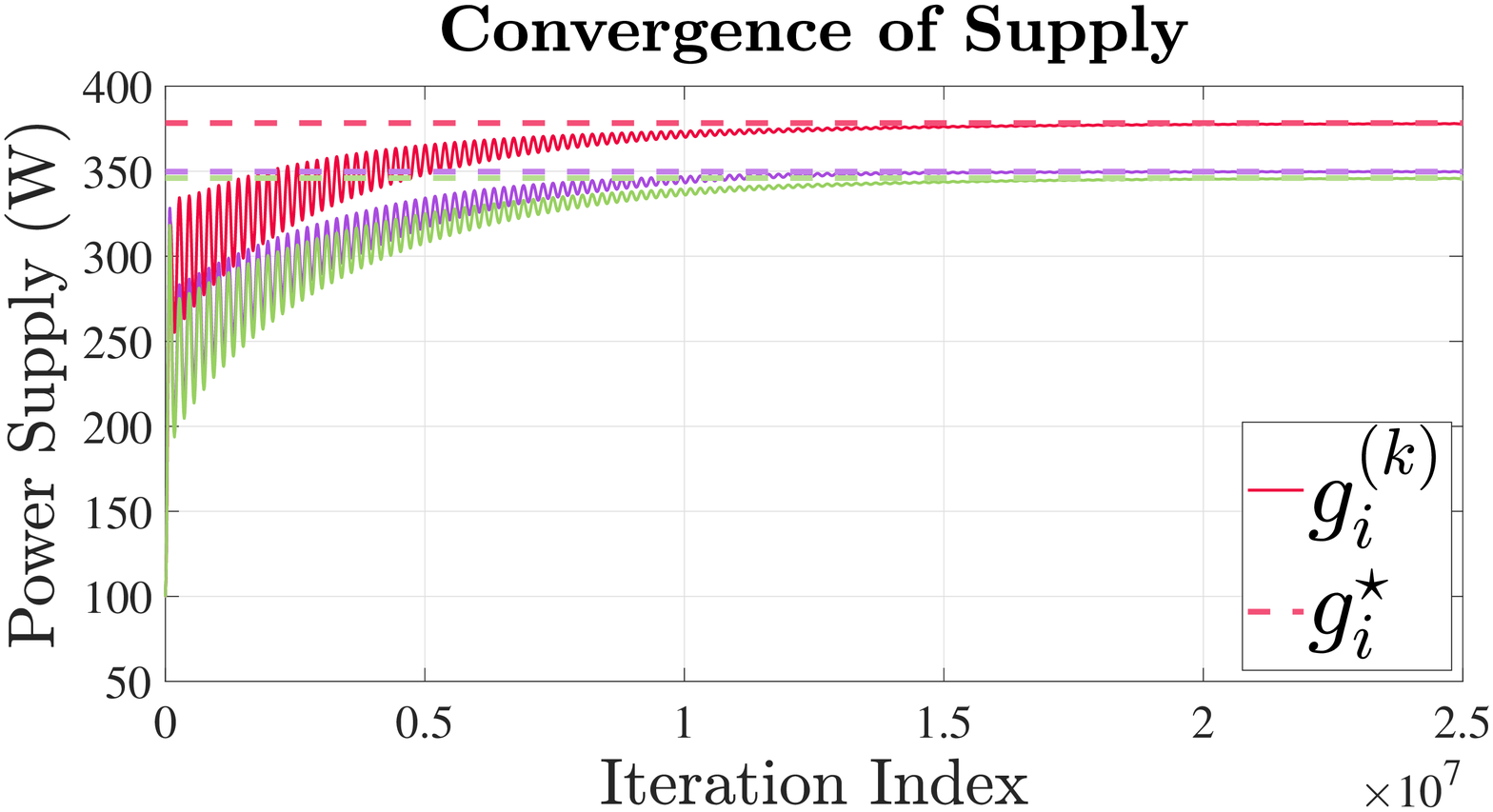}
         \caption{}
         \label{fig:powernetworksupply}
     \end{subfigure}
        \caption{Numerical study results for power network under dynamic impersonation attack: (a)/(b) displaying convergence of the demand/supply, respectively.}
        \label{fig:powernetresults}
\end{figure}

Observe that the formulation in \eqref{eq:poweroptimization} does not directly match with our general formulation in \eqref{eq:mainopt}  mainly due to the presence of equality constraint \eqref{const:supplydemand}, which prevents the application of the robustified formulation in \eqref{eq:robust} and hence the robust PD-DRA algorithm for static impersonation attacks. Nevertheless, our algorithm for dynamic impersonation attacks can still be applied since it does not require any robustified constraints (which cannot be done for equality constraints).


We have chosen the utility function for load $i$ to be $U_i(d_i)=\beta_i \log d_i$ and randomly generated $\beta_i$ from a uniform distribution in $[500,1000]$. For generators, we set the cost function $C_i(g_i)=e^{c_i g_i}$, where $c_1=0.01,~c_2=0.011,~c_3=0.012$. We obtained the Power Transfer Distribution Factor (PTDF) matrix \textbf{H} and the vector of flow limits $\boldsymbol{c}$ from \textsc{Matpower}\cite{matpower}. To simulate a dynamic impersonation attack scenario, we assigned a probability $p$ for an uplink to be compromised at each iteration. We ran Algorithm~\ref{alg:jamming_alg} for $p=0.15$, $m=75$ and $\alpha_2=0.49$.

The results are shown in Figure~\ref{fig:powernetresults}. \bt{In both Figures~\ref{fig:powernetworkdemand} and \ref{fig:powernetworksupply}, each different color corresponds to a different agent. A colored solid line corresponds to an agent's parameter iterates generated by the algorithm. A dashed line with the same color as a solid line is the to the optimal value of that agent's parameter, which is the solution of the regularized optimization problem (formulated as \eqref{eq:pdreg}) of \eqref{eq:poweroptimization}.} Our algorithm successfully generates sequences that convergence to the optimal solution of the regularized problem for both power supplying and power demanding agents.



\section{Conclusion}
In this paper, we studied two strategies for establishing primal-dual algorithms for resource allocation in presence of Byzantine attackers. Specifically, we consider static and dynamic impersonation attack scenarios and propose an attack-resilient primal-dual algorithm for each scenario based on robust mean estimation techniques. We derive bounds for the performance (in terms of distance to optimality) of the proposed algorithms and show that our algorithm for static impersonation attack converges to a neighborhood of the optimal solution of the regularized and robustified resource allocation problem, whereas our algorithm for dynamic impersonation attack converges to the optimal solution of the original regularized problem. We verify our theoretical results via computational simulations for network utility maximization problems involving optimal distributed resource allocation, such as power distribution networks.
\bibliographystyle{IEEEtran}
\bibliography{references,robust_bib,num_bib,opt_bib}

\begin{thebibliography}{10}
\providecommand{\url}[1]{#1}
\csname url@samestyle\endcsname
\providecommand{\newblock}{\relax}
\providecommand{\bibinfo}[2]{#2}
\providecommand{\BIBentrySTDinterwordspacing}{\spaceskip=0pt\relax}
\providecommand{\BIBentryALTinterwordstretchfactor}{4}
\providecommand{\BIBentryALTinterwordspacing}{\spaceskip=\fontdimen2\font plus
\BIBentryALTinterwordstretchfactor\fontdimen3\font minus
  \fontdimen4\font\relax}
\providecommand{\BIBforeignlanguage}[2]{{%
\expandafter\ifx\csname l@#1\endcsname\relax
\typeout{** WARNING: IEEEtran.bst: No hyphenation pattern has been}%
\typeout{** loaded for the language `#1'. Using the pattern for}%
\typeout{** the default language instead.}%
\else
\language=\csname l@#1\endcsname
\fi
#2}}
\providecommand{\BIBdecl}{\relax}
\BIBdecl

\bibitem{kelly}
F.~P. Kelly, A.~K. Maulloo, and D.~K.~H. Tan, ``Rate control for communication
  networks: shadow prices, proportional fairness and stability,'' \emph{J OPER
  RES SOC}, vol.~49, no.~3, pp. 237--252, Mar 1998.

\bibitem{low}
S.~H. {Low} and D.~E. {Lapsley}, ``Optimization flow control. i. basic
  algorithm and convergence,'' \emph{IEEE/ACM Transactions on Networking},
  vol.~7, no.~6, pp. 861--874, Dec 1999.

\bibitem{mohsenian}
P.~{Samadi}, A.~{Mohsenian-Rad}, R.~{Schober}, V.~W.~S. {Wong}, and
  J.~{Jatskevich}, ``Optimal real-time pricing algorithm based on utility
  maximization for smart grid,'' in \emph{IEEE SmartGridComm}, Oct 2010, pp.
  415--420.

\bibitem{lina}
N.~{Li}, L.~{Chen}, and S.~H. {Low}, ``Optimal demand response based on utility
  maximization in power networks,'' in \emph{2011 IEEE Power and Energy Society
  General Meeting}, July 2011, pp. 1--8.

\bibitem{mung}
M.~Chiang and J.~{Bell}, ``Balancing supply and demand of bandwidth in wireless
  cellular networks: utility maximization over powers and rates,'' in
  \emph{IEEE INFOCOM 2004}, vol.~4, March 2004, pp. 2800--2811 vol.4.

\bibitem{cache}
M.~{Dehghan}, L.~{Massoulie}, D.~{Towsley}, D.~{Menasche}, and Y.~C. {Tay}, ``A
  utility optimization approach to network cache design,'' in \emph{IEEE
  INFOCOM}, April 2016, pp. 1--9.

\bibitem{wsn1}
M.~Zhao, J.~Li, and Y.~Yang, ``Joint mobile energy replenishment and data
  gathering in wireless rechargeable sensor networks,'' in \emph{Proceedings of
  the 23rd International Teletraffic Congress}, ser. ITC '11.\hskip 1em plus
  0.5em minus 0.4em\relax International Teletraffic Congress, 2011, pp.
  238--245.

\bibitem{wsn2}
R.~{Deng}, Y.~{Zhang}, S.~{He}, J.~{Chen}, and X.~{Shen}, ``Maximizing network
  utility of rechargeable sensor networks with spatiotemporally coupled
  constraints,'' \emph{IEEE JSAC}, vol.~34, no.~5, pp. 1307--1319, May 2016.

\bibitem{traffic}
N.~{Mehr}, J.~{Lioris}, R.~{Horowitz}, and R.~{Pedarsani}, ``Joint perimeter
  and signal control of urban traffic via network utility maximization,'' in
  \emph{IEEE ITSC}, Oct 2017, pp. 1--6.

\bibitem{koshal2011multiuser}
J.~Koshal, A.~Nedi{\'c}, and U.~V. Shanbhag, ``Multiuser optimization:
  Distributed algorithms and error analysis,'' \emph{SIAM Journal on
  Optimization}, vol.~21, no.~3, pp. 1046--1081, 2011.

\bibitem{palomar2006tutorial}
D.~P. Palomar and M.~Chiang, ``A tutorial on decomposition methods for network
  utility maximization,'' \emph{IEEE JSAC}, vol.~24, no.~8, pp. 1439--1451,
  2006.

\bibitem{sundaram2011distributed}
S.~Sundaram and C.~N. Hadjicostis, ``Distributed function calculation via
  linear iterative strategies in the presence of malicious agents,'' \emph{IEEE
  TAC}, vol.~56, no.~7, pp. 1495--1508, 2011.

\bibitem{pasqualetti2012consensus}
F.~Pasqualetti, A.~Bicchi, and F.~Bullo, ``Consensus computation in unreliable
  networks: A system theoretic approach,'' \emph{IEEE TAC}, vol.~57, no.~1, pp.
  90--104, 2012.

\bibitem{gentz2016data}
R.~Gentz, S.~X. Wu, H.-T. Wai, A.~Scaglione, and A.~Leshem, ``Data injection
  attacks in randomized gossiping,'' \emph{IEEE TSIPN}, vol.~2, no.~4, pp.
  523--538, 2016.

\bibitem{sundaram2018distributed}
S.~Sundaram and B.~Gharesifard, ``Distributed optimization under adversarial
  nodes,'' \emph{IEEE TAC}, 2018.

\bibitem{chen2018resilient}
Y.~Chen, S.~Kar, and J.~Moura, ``Resilient distributed estimation: Sensor
  attacks,'' \emph{IEEE Transactions on Automatic Control}, 2018.

\bibitem{consensusbenameur}
W.~Ben-Ameur, P.~Bianchi, and J.~Jakubowicz, ``Robust distributed consensus
  using total variation,'' \emph{IEEE Transactions on Automatic Control},
  vol.~61, no.~6, pp. 1550--1564, 2016.

\bibitem{consensusleblanc}
H.~J. {LeBlanc}, H.~{Zhang}, X.~{Koutsoukos}, and S.~{Sundaram}, ``Resilient
  asymptotic consensus in robust networks,'' \emph{IEEE Journal on Selected
  Areas in Communications}, vol.~31, no.~4, pp. 766--781, April 2013.

\bibitem{consensusbaras}
J.~S. {Baras} and X.~{Liu}, ``Trust is the cure to distributed consensus with
  adversaries,'' in \emph{2019 27th Mediterranean Conference on Control and
  Automation (MED)}, July 2019, pp. 195--202.

\bibitem{ravi19consensus}
N.~{Ravi}, A.~{Scaglione}, and A.~{Nedić}, ``A case of distributed
  optimization in adversarial environment,'' in \emph{2019 IEEE International
  Conference on Acoustics, Speech and Signal Processing}, May 2019, pp.
  5252--5256.

\bibitem{feng2014distributed}
J.~Feng, H.~Xu, and S.~Mannor, ``Distributed robust learning,'' \emph{arXiv
  preprint arXiv:1409.5937}, 2014.

\bibitem{yin2018byzantine}
D.~Yin, Y.~Chen, K.~Ramchandran, and P.~Bartlett, ``Byzantine-robust
  distributed learning: Towards optimal statistical rates,'' \emph{arXiv
  preprint arXiv:1803.01498}, 2018.

\bibitem{alistarh2018byzantine}
D.~Alistarh, Z.~Allen-Zhu, and J.~Li, ``Byzantine stochastic gradient
  descent,'' in \emph{NeurIPS}, 2018, pp. 4618--4628.

\bibitem{chen2017ml}
Y.~Chen, L.~Su, and J.~Xu, ``Distributed statistical machine learning in
  adversarial settings: Byzantine gradient descent,'' \emph{Proc. ACM Meas.
  Anal. Comput. Syst.}, vol.~1, no.~2, pp. 44:1--44:25, 2017.

\bibitem{data19ml}
D.~{Data}, L.~{Song}, and S.~{Diggavi}, ``Data encoding methods for
  byzantine-resilient distributed optimization,'' in \emph{2019 IEEE
  International Symposium on Information Theory (ISIT)}, July 2019, pp.
  2719--2723.

\bibitem{Donoho1983}
D.~L. Donoho and P.~J. Huber, ``The notion of breakdown point,'' \emph{A
  festschrift for Erich L. Lehmann}, vol. 157184, 1983.

\bibitem{huber2011robust}
P.~J. Huber, \emph{Robust statistics}.\hskip 1em plus 0.5em minus 0.4em\relax
  Springer, 2011.

\bibitem{minsker2015geometric}
S.~Minsker \emph{et~al.}, ``Geometric median and robust estimation in banach
  spaces,'' \emph{Bernoulli}, vol.~21, no.~4, pp. 2308--2335, 2015.

\bibitem{diakonikolas2016robust}
I.~Diakonikolas, G.~Kamath, D.~M. Kane, J.~Li, A.~Moitra, and A.~Stewart,
  ``Robust estimators in high dimensions without the computational
  intractability,'' in \emph{IEEE FOCS}, 2016, pp. 655--664.

\bibitem{Steinhardt2017}
J.~Steinhardt, M.~Charikar, and G.~Valiant, ``Resilience: A criterion for
  learning in the presence of arbitrary outliers,'' \emph{arXiv preprint
  arXiv:1703.04940}, 2017.

\bibitem{cdcversion}
C.~A. {Uribe}, H.~{Wai}, and M.~{Alizadeh}, ``Resilient distributed
  optimization algorithms for resource allocation,'' in \emph{2019 IEEE 58th
  Conference on Decision and Control (CDC)}, 2019, pp. 8341--8346.

\bibitem{nesterov2005smooth}
Y.~Nesterov, ``Smooth minimization of non-smooth functions,''
  \emph{Mathematical programming}, vol. 103, no.~1, pp. 127--152, 2005.

\bibitem{uribe2018dual}
C.~A. Uribe, S.~Lee, A.~Gasnikov, and A.~Nedi{\'c}, ``A dual approach for
  optimal algorithms in distributed optimization over networks,'' \emph{arXiv
  preprint arXiv:1809.00710}, 2018.

\bibitem{onlineversion}
B.~Turan, C.~A. Uribe, H.-T. Wai, and M.~Alizadeh, ``Resilient primal-dual
  optimization algorithms for distributed resource allocation,'' \emph{arXiv
  preprint arXiv:2001.00612}, 2020.

\bibitem{iag_blatt}
D.~Blatt, A.~O. Hero, and H.~Gauchman, ``A convergent incremental gradient
  method with a constant step size,'' \emph{SIAM Journal on Optimization},
  vol.~18, no.~1, pp. 29--51, 2007.

\bibitem{iag_metasupablo}
M.~Gürbüzbalaban, A.~Ozdaglar, and P.~A. Parrilo, ``On the convergence rate
  of incremental aggregated gradient algorithms,'' \emph{SIAM Journal on
  Optimization}, vol.~27, no.~2, pp. 1035--1048, 2017.

\bibitem{iag_Tseng2014}
P.~Tseng and S.~Yun, ``Incrementally updated gradient methods for constrained
  and regularized optimization,'' \emph{Journal of Optimization Theory and
  Applications}, vol. 160, no.~3, pp. 832--853, Mar 2014.

\bibitem{cvx}
M.~Grant and S.~Boyd, ``{CVX}: Matlab software for disciplined convex
  programming, version 2.1,'' \url{http://cvxr.com/cvx}, Mar. 2014.

\bibitem{matpower}
R.~D. {Zimmerman}, C.~E. {Murillo-Sánchez}, and R.~J. {Thomas}, ``Matpower:
  Steady-state operations, planning, and analysis tools for power systems
  research and education,'' \emph{IEEE Transactions on Power Systems}, vol.~26,
  no.~1, pp. 12--19, Feb 2011.

\bibitem{mjohansson}
H.~R. {Feyzmahdavian}, A.~{Aytekin}, and M.~{Johansson}, ``A delayed proximal
  gradient method with linear convergence rate,'' in \emph{2014 IEEE
  International Workshop on Machine Learning for Signal Processing (MLSP)},
  Sep. 2014, pp. 1--6.

\end{thebibliography}
\newpage
\begin{IEEEbiography}[{\includegraphics[width=1in,height=1.21in,clip,keepaspectratio]{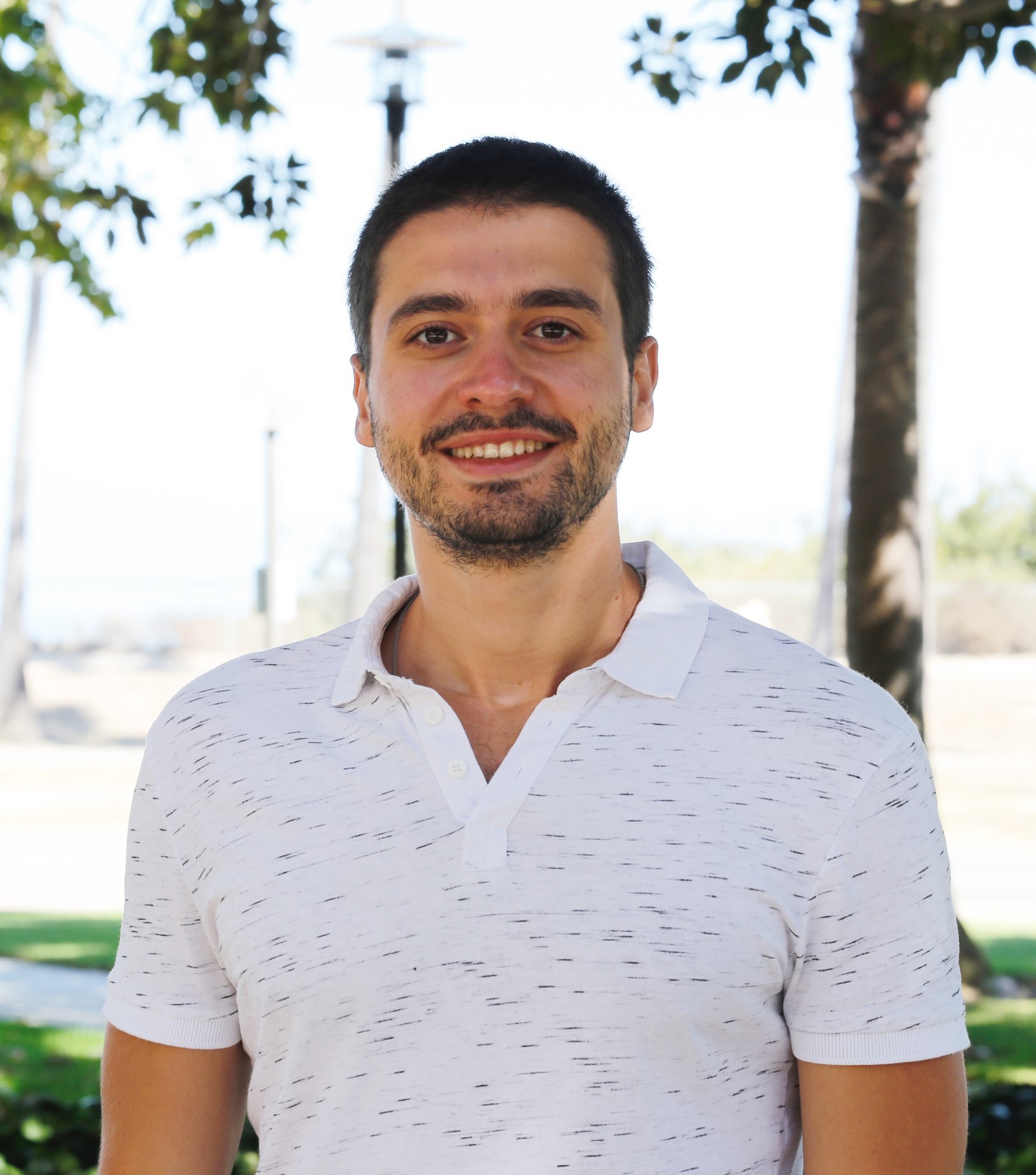}}]{BERKAY TURAN}is pursuing the Ph.D. degree in Electrical and Computer Engineering at the University of California, Santa Barbara. He received the B.Sc. degree in Electrical and Electronics Engineering as well as the B.Sc. degree in  Physics degree from \ Bo\u gazi\c ci University, Istanbul, Turkey, in 2018. His research interests include optimization and reinforcement learning for the design, control, and analysis of smart infrastructure systems such as the power grid and transportation systems.
\end{IEEEbiography}
\vskip -2\baselineskip plus -1fil
\begin{IEEEbiography}[{\includegraphics[width=1in,height=1.21in,clip,keepaspectratio]{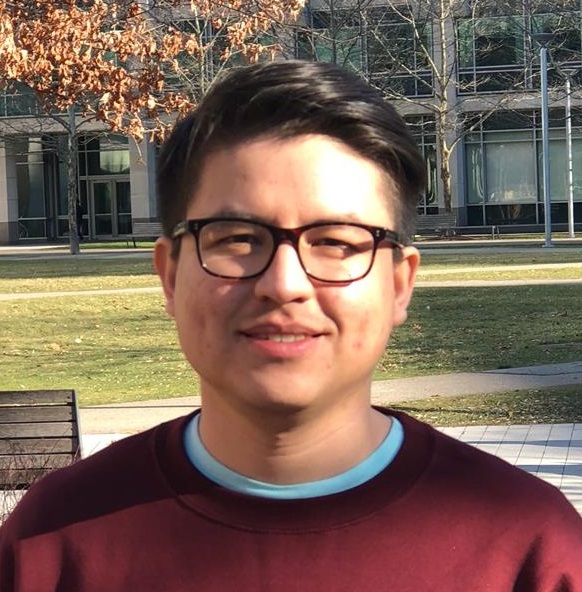}}]{C\'ESAR A. URIBE} received the M.Sc. degrees in systems and control from Delft University of Technology, in The Netherlands, and in applied mathematics from the University of Illinois at Urbana-Champaign, in 2013 and 2016, respectively. He also received the PhD degree in electrical and computer engineering at the University of Illinois at Urbana-Champaign in 2018.  He is currently a Postdoctoral Associate in the Laboratory for Information and Decision Systems-LIDS at the Massachusetts Institute of Technology-MIT. His research interests include distributed learning and optimization, decentralized control, algorithm analysis, and computational optimal transport.
\end{IEEEbiography}
\vskip -2\baselineskip plus -1fil
\begin{IEEEbiography}[{\includegraphics[width=1in,height=1.21in,clip,keepaspectratio]{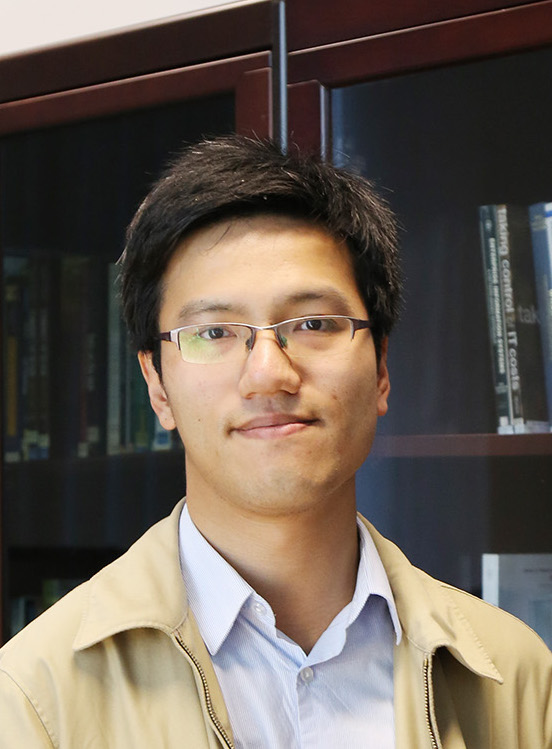}}]{Hoi-To Wai} (S’11–M’18) received his PhD degree from Arizona State University (ASU) in Electrical Engineering in Fall 2017, B. Eng. (with First Class Honor) and M. Phil. degrees in Electronic Engineering from The Chinese University of Hong Kong (CUHK) in 2010 and 2012, respectively. He is an Assistant Professor in the Department of Systems Engineering \& Engineering Management at CUHK. He has held research positions at ASU, UC Davis, Telecom ParisTech, Ecole Polytechnique, LIDS, MIT.
Hoi-To's research interests are in the broad area of signal processing, machine learning and distributed optimization, with a focus of their applications to network science. His dissertation has received the 2017's Dean's Dissertation Award from the Ira A. Fulton Schools of Engineering of ASU and he is a recipient of a Best Student Paper Award at ICASSP 2018.
\end{IEEEbiography}
\vskip -2\baselineskip plus -1fil
\begin{IEEEbiography}[{\includegraphics[width=1in,height=1.21in,clip,keepaspectratio]{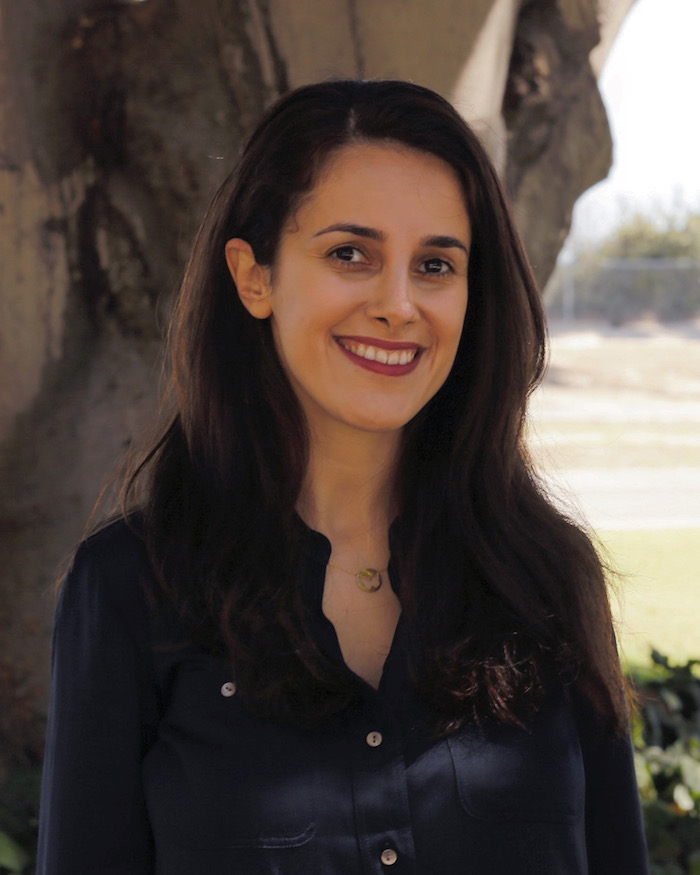}}]{MAHNOOSH ALIZADEH} is an assistant professor of Electrical and Computer Engineering at the University of California Santa Barbara. She received the B.Sc. degree (’09) in Electrical Engineering from Sharif University of Technology and the M.Sc. (’13) and Ph.D. (’14) degrees in Electrical and Computer Engineering from the University of California Davis, where she was the recipient of the Richard C. Dorf award for outstanding research accomplishment. From 2014 to 2016, she was a postdoctoral scholar at Stanford University. Her research interests are focused on designing network control, optimization, and learning frameworks to promote efficiency and resiliency in societal-scale cyber-physical systems. Dr. Alizadeh is a recipient of the NSF CAREER award and the best paper award from HICSS-53 power systems track.
\end{IEEEbiography}
\newpage
\appendix
\subsection{Proof of Lemma \ref{lem:conservativeform}}\label{app:conservative}
\bt{Observe that in both \eqref{eq:robust} and \eqref{eq:reformulate_s}, the decision variables of the optimization problems are $\{\prm_i\}_{i\in{\cal H}}$. Hence, it suffices to show that any given set of $\{\prm_i\}_{i\in{\cal H}}$ satisfying constraints of \eqref{eq:reformulate_s} also satisfies constraints of \eqref{eq:robust}.}
Let $\overline{\prm}_\mathcal{A}\eqdef\frac{1}{|\mathcal{A}|}\sum_{i\in\mathcal A}\prm_i$.
Since $g_t$ is $L$-smooth, the following holds:
\bt{
\beq
\begin{split}\label{eq:robustconstraintproof}
 &\underset{\prm_j \in \Cset_j,j\in \mathcal A}{\max} \textstyle g_t \big( \frac{1}{N} \sum_{i=1}^N \prm_i \big)=\underset{\prm_j \in \Cset_j,j\in \mathcal A}{\max} g_t\big(\frac{|\mathcal H|}{N}\overline{\prm}_\mathcal{H}+\frac{|\mathcal A|}{N}\overline{\prm}_\mathcal{A})    \\
 &\textstyle =\underset{\prm_j \in \Cset_j,j\in \mathcal A}{\max}  g_t \big( (1-\alpha_1)\overline{\prm}_\mathcal{H}+(\alpha_1-\frac{|\mathcal A|}{N})\overline{\prm}_\mathcal{H}+\frac{|\mathcal A|}{N}\overline{\prm}_\mathcal{A}\big)\\
 &\textstyle \leq\underset{\prm_j \in \Cset_j,j\in \mathcal A}{\max}  \bigg(g_t \big( (1-\alpha_1)\overline{\prm}_\mathcal{H}\big)\\
 &\hspace{2cm}+\Big\langle \widetilde{\prm}, \grd g_t \big((1-\alpha_1) \overline{\prm}_\mathcal{H}\big) \Big\rangle +\frac{L}{2}\|\widetilde{\prm}\|^2\bigg),
\end{split}
\eeq}
where we defined $\widetilde{\prm}\eqdef(\alpha_1-\frac{|\mathcal A|}{N})\overline{\prm}_\mathcal{H}+\frac{|\mathcal A|}{N}\overline{\prm}_\mathcal{A}$. Observe that
\begin{equation}
    \textstyle \|\widetilde{\prm}\|\leq (\alpha_1-\frac{|\mathcal A|}{N})\|\overline{\prm}_\mathcal{H}\|+\frac{|\mathcal A|}{N}\|\overline{\prm}_\mathcal{A}\|\leq \alpha_1 R.
\end{equation}
Furthermore, since the gradient of $g_t$ is uniformly bounded by $B$ and $\alpha_1^2\leq\alpha_1$, \eqref{eq:robustconstraintproof} can be upper bounded by:
\beq
g_t \big((1-\alpha_1)\overline{\prm}_\mathcal{H}\big) + \alpha_1 \big( R B + {\textstyle \frac{1}{2}} L R^2 \big)
\eeq
\bt{Hence, we have shown that for any given set of $\{\prm_i\}_{i\in{\cal H}}$, the following holds:
\beq
\underset{\prm_j \in \Cset_j,j\in \mathcal A}{\max} g_t \big( {\textstyle \frac{1}{N} \sum_{i=1}^N \prm_i } \big) \leq g_t \big((1-\alpha_1)\overline{\prm}_\mathcal{H}\big) + \alpha_1 \big( R B + {\textstyle \frac{1}{2}} L R^2 \big).
\eeq}
As such defining $c_t \eqdef \alpha_1\big( R B + \frac{1}{2} L R^2 \big)$,
it can be seen that \bt{if a set of $\{\prm_i\}_{i\in{\cal H}}$ satisfies
\beq
g_t \big( (1-\alpha_1)\overline{\prm}_\mathcal{H} \big) + c_t \leq 0,~t=1,...,T,
\eeq
    the same set of $\{\prm_i\}_{i\in{\cal H}}$ satisfies the desired constraint \eqref{eq:robust}.} 
\subsection{Proof of Proposition \ref{prop:median}}\label{app:medianest}
Fix any $j \in [d]$. The assumption implies that for all $i \in {\cal H}$, one has:
\begin{equation}\label{eq:boundinftynorm}
    | [ \bm x_i - \overline{\bm x}_{\cal H} ]_j| \leq r.
\end{equation}
We observe that $|{\cal H}| \geq (1-\alpha_1) N$. Applying \cite[Lemma 1]{feng2014distributed} shows that the median estimator\footnote{At each coordinate, the median is the geometric median estimator of one dimension in \cite{feng2014distributed}.} satisfies  
\begin{equation}
| [ {\bm x}_{\sf med} -  \overline{\bm x}_{\cal H} ]_j | \leq (1- \alpha_1) \sqrt{\frac{1}{1-2\alpha_1}} r.
\end{equation}
The above implies that for all $i \in {\cal H}$, we have
\begin{equation}
| [ \bm x_i - {\bm x}_{\sf med} ]_j | \leq \left( 1 + \sqrt{\frac{(1-\alpha_1)^2}{1-2\alpha_1}} \right) r.    
\end{equation}
This implies that $r_j \leq \left( 1 + \sqrt{\frac{(1-\alpha_1)^2}{1-2\alpha_1}} \right) r$, since $|{\cal H}| \geq (1-\alpha_1)N$.
We then bound the performance of $\widehat{\bm x}_{\mathcal H}$:
\begin{equation}
\begin{split}
& (1-\alpha_1)N [ \widehat{\bm x}_\mathcal H ]_j = \sum_{i \in {\cal N}_j}  [{\bm x}_i ]_j \\
& =
\sum_{i \in {\cal H} }  [{\bm x}_i ]_j - \hspace{-.2cm} \sum_{i \in {\cal H} \setminus {\cal N}_j }  [{\bm x}_i ]_j + \hspace{-.2cm} \sum_{i \in {\cal A} \cap {\cal N}_j }  [{\bm x}_i ]_j,
\end{split}    
\end{equation}
thus
\begin{equation}
\begin{split}
& (1-\alpha_1)N [ \widehat{\bm x}_{\cal H}- \overline{\bm x}_{\cal H} ]_j  \\
& = - \sum_{i \in {\cal H} \setminus {\cal N}_j }  [{\bm x}_i - \overline{\bm x }_{\cal H} ]_j + \hspace{-.2cm} \sum_{i \in {\cal A} \cap {\cal N}_j }  [{\bm x}_i - \overline{\bm x}_{\cal H} ]_j\\
& = - \sum_{i \in {\cal H} \setminus {\cal N}_j }  [{\bm x}_i - \overline{\bm x }_{\cal H} ]_j + \hspace{-.2cm} \sum_{i \in {\cal A} \cap {\cal N}_j }  [{\bm x}_i - {\bm x}_{\sf med} + {\bm x}_{\sf med} - \overline{\bm x}_{\cal H} ]_j
\end{split}    
\end{equation}
Notice that $|{\cal A} \cap {\cal N}_j| \leq \alpha_1 N$ and thus $|{\cal H} \setminus {\cal N}_j| \leq \alpha_1 N$. Gathering terms shows
\begin{equation}
| [ \widehat{\bm x} - \overline{\bm x}_{\cal H} ]_j | \leq \frac{2\alpha_1 N}{(1-\alpha_1)N}
\left( 1 + \sqrt{\frac{(1-\alpha_1)^2}{1-2\alpha_1}} \right) r.    
\end{equation}

The above holds for all $j \in [d]$. Applying the norm equivalence shows the desired bound.
\subsection{Proof of Lemma \ref{lem:perturbedgradients1}}\label{app:permanentperturbedgrads}
\bt{
Let $\big[ {\bm e}_{\prm}^{(k)} \big]_i$ denote the $i$th block of ${\bm e}_{\prm}^{(k)}$, and $\big[\widehat{\bm g}_{ \prm }^{(k)}\big]_i$ denote the $i$th block of $\widehat{\bm g}_{ \prm }^{(k)}$. From Equation \eqref{eq:primal_pb}: 
\beq
\big[\widehat{\bm g}_{ \prm }^{(k)}\big]_i=\frac{1}{N} \big( 
\widehat{\bm g}^{(k)}_{\cal H} + \grd_{\prm_i} f_i( \prm_i^{(k)} ) + \upsilon \prm_i^{(k)} \big).
\eeq
Furthermore, we replace $\widehat{\bm g}^{(k)}_{\cal H}$ from Algorithm~\ref{alg:robust_pddra} Step~\ref{r_step:accum0}\ref{r_step:broadcast}: 
\beq
\label{eq:prm_perturbed_grad}
\big[\widehat{\bm g}_{ \prm }^{(k)}\big]_i=\frac{1}{N} \Big( 
\sum_{t=1}^T \lambda_t^{(k)} \grd_{\prm} \overline{g}_t ((1-\alpha_1)\widehat{\prm}_{\cal H}^{(k)})  + \grd_{\prm_i} f_i( \prm_i^{(k)} ) + \upsilon \prm_i^{(k)} \Big).
\eeq
The perturbation $\big[{\bm e}_{\prm}^{(k)}\big]_i=\big[\widehat{\bm g}_{ \prm }^{(k)}\big]_i-\grd_{ \prm_i } \overline{\cal L}_{\upsilon} ( \{ \prm_i \}_{i \in {\cal H}}; \bm{\lambda}; {\cal H} )$ is given by the difference between \eqref{eq:prm_perturbed_grad} and \eqref{eq:grd_prm_robust}:
\beq
\begin{split}
    \big[{\bm e}_{\prm}^{(k)}\big]_i=\frac{1}{N}\sum_{t=1}^T &\lambda_t^{(k)}\Big( \grd_{\prm} \overline{g}_t ((1-\alpha_1)\widehat{\prm}_{\cal H}^{(k)})\\
    &-\frac{(1-\alpha_1)N}{|\mathcal H|}\grd_{\prm} \overline{g}_t ((1-\alpha_1)\bar{\prm}^{(k)}_\mathcal{H})\Big).
\end{split}
\eeq
By adding and subtracting $\frac{1}{N}\Big(\sum_{t=1}^T \lambda_t^{(k)}\grd_{\prm} \overline{g}_t ((1-\alpha_1)\bar{\prm}^{(k)}_\mathcal{H}\Big)$, the above expression becomes:
\beq \label{eq:etheta}
\begin{aligned}
\big[ {\bm e}_{\prm}^{(k)} \big]_i = \frac{1}{N}  {\sum_{t=1}^T} \lambda_t^{(k)} \Big(& \grd_{\prm} \overline{g}_t \big( {\textstyle (1-\alpha_1) } \widehat{\prm}_{\cal H}^{(k)} \big) \\&\hspace{-1.5cm}-  \grd_{\prm} \overline{g}_t \big( {\textstyle (1-\alpha_1) } \overline{\prm}_{\cal H}^{(k)} \big)\\
&\hspace{-1.5cm}+\frac{|\mathcal H|-(1-\alpha_1)N}{|\mathcal H|}\grd_{\prm} \overline{g}_t \big( {\textstyle (1-\alpha_1) } \overline{\prm}_{\cal H}^{(k)}\big)\Big) 
\end{aligned}
\eeq}

Similarly, comparing \eqref{eq:per_lam} with \eqref{eq:dual_pb} and \eqref{eq:grd_lam_robust}, we identify that:

\beq \label{eq:elambda}
\big[ {\bm e}_{\bm{\lambda}}^{(k)} \big]_t = \overline{g}_t \big( {\textstyle  (1-\alpha_1) }\widehat{\prm}_{\cal H}^{(k)} \big) - \overline{g}_t \big( {\textstyle  (1-\alpha_1) }\overline{\prm}_{\cal H}^{(k)} \big).
\eeq

Using Assumption~\ref{ass:bdd} and the said assumptions, we immediately see that
\beq
\begin{aligned}
\| \big[ {\bm e}_{\prm}^{(k)} \big]_i \| \leq& (1-\alpha_1) \frac{\overline{\lambda}L T}{N}  \| 
\widehat{\prm}_{\cal H}^{(k)} - \overline{\prm}_{\cal H}^{(k)} \|\\
&+\frac{|\mathcal H|-(1-\alpha_1)N}{|\mathcal H|}\frac{\overline{\lambda}B T}{N}
\end{aligned}
\eeq
which then implies \eqref{eq:bd_prm}.
 Assumption~\ref{ass:bdd} implies that $\overline{g}_t$ is $B$-Lipschitz continuous,
therefore 
\beq
| \big[ {\bm e}_{\bm{\lambda}}^{(k)} \big]_t | \leq B(1-\alpha_1) \| 
\widehat{\prm}_{\cal H}^{(k)} - \overline{\prm}_{\cal H}^{(k)} \|,
\eeq
which implies \eqref{eq:bd_lambda}.
\subsection{Proof of Theorem \ref{thm:main}}\label{app:permanentthm}
Based on Lemma~\ref{lem:perturbedgradients1}, our idea is to perform a perturbation analysis on the PDA algorithm.
Without loss of generality, we assume $N=1$ and denote $\prm = \prm_1$. 
To simplify notations we define $\upsilon'\eqdef(1-\alpha_1)\upsilon$. We also drop the subscript, denote the modified and regularized Lagrangian function as ${\cal L} = \overline{\cal L}_{\upsilon}$.
Furthermore, we denote the saddle point to \eqref{eq:pdp} as ${\bm z}^\star = (\prm^\star, \bm{\lambda}^\star)$. 

Using the fact that $\prm^\star = {\cal P}_{\Cset}( \prm^\star  )= {\cal P}_{\Cset} \big( \prm^\star - \gamma \grd_{\prm} {\cal L} ( \prm^\star, \bm{\lambda}^\star) \big)$, we 
observe that in the primal update:
\begin{align}
\nonumber& \| \prm^{(k+1)} - \prm^\star \|^2 \label{eq:ss} \\
\nonumber& \overset{(a)}{\leq} \| \prm^{(k)} - \prm^\star \|^2 - 2 \gamma 
\langle \widehat{\bm g}_{\prm}^{(k)} - \grd_{\prm} {\cal L} ( \prm^\star, \bm{\lambda}^\star), \prm^{(k)} - \prm^\star \rangle \\
& \hspace{.8cm} + \gamma^2 
\| \widehat{\bm g}_{\prm}^{(k)} - \grd_{\prm} {\cal L} ( \prm^\star, \bm{\lambda}^\star) \|^2 
\end{align}
where (a)  is due to the projection inequality $\| {\cal P}_{\Cset} ( {\bm x} - {\bm y} ) \| \leq 
\| {\bm x} - {\bm y} \|$. Furthermore, using the Young's inequality, for any $c_0, c_1 > 0$, we have
\begin{align}
\nonumber& \| \prm^{(k+1)} - \prm^\star \|^2 \\
\nonumber& \leq \| \prm^{(k)} - \prm^\star \|^2 \\
\nonumber& \hspace{.5cm} - 2 \gamma \langle \grd_{\prm} {\cal L}( \prm^{(k)}, \bm{\lambda}^{(k)} ) - \grd_{\prm} {\cal L} ( \prm^\star, \bm{\lambda}^\star), \prm^{(k)} - \prm^\star \rangle \\
\nonumber& \hspace{.5cm} { + \gamma^2 ( 1 + c_0 ) \|  \grd_{\prm} {\cal L}( \prm^{(k)}, \bm{\lambda}^{(k)} ) - \grd_{\prm} {\cal L} ( \prm^\star, \bm{\lambda}^\star) \|^2} \\
\nonumber& \hspace{.5cm} {- 2 \gamma \langle {\bm e}_{\prm}^{(k)},  \prm^{(k)} - \prm^\star \rangle + \gamma^2 \big( 1 + \frac{1}{c_0 } \big) \| {\bm e}_{\prm}^{(k)} \|^2} \\
\nonumber& \leq (1 + 2 c_1 \gamma ) \!~  \| \prm^{(k)} - \prm^\star \|^2 \\
\nonumber& \hspace{.5cm} - 2 \gamma \langle \grd_{\prm} {\cal L}( \prm^{(k)}, \bm{\lambda}^{(k)} ) - \grd_{\prm} {\cal L} ( \prm^\star, \bm{\lambda}^\star), \prm^{(k)} - \prm^\star \rangle \\
\nonumber& \hspace{.5cm} + \gamma^2 ( 1 + c_0 ) \|  \grd_{\prm} {\cal L}( \prm^{(k)}, \bm{\lambda}^{(k)} ) - \grd_{\prm} {\cal L} ( \prm^\star, \bm{\lambda}^\star) \|^2  \\
& \hspace{.5cm} + \Big( \frac{2 \gamma}{c_1} + \gamma^2 + \frac{\gamma^2}{c_0 } \Big) \| {\bm e}_{\prm}^{(k)} \|^2 .
\end{align}
Similarly, in the dual update we get, 
\begin{align}
\nonumber& \| \bm{\lambda}^{(k+1)} - \bm{\lambda}^\star \|^2 \\
\nonumber& \leq \| \bm{\lambda}^{(k)} - \bm{\lambda}^\star \|^2 + \gamma^2 
\| \widehat{\bm g}_{\bm{\lambda}}^{(k)} - \grd_{\bm{\lambda}} {\cal L} ( \prm^\star, \bm{\lambda}^\star) \|^2\\
\nonumber& \hspace{.5cm} + 2 \gamma 
\langle \widehat{\bm g}_{\bm{\lambda}}^{(k)} - \grd_{\bm{\lambda}} {\cal L} ( \prm^\star, \bm{\lambda}^\star), \bm{\lambda}^{(k)} - \bm{\lambda}^\star \rangle \\
\nonumber& \leq (1 + 2 c_1 \gamma ) \!~\| \bm{\lambda}^{(k)} - \bm{\lambda}^\star \|^2 \\
\nonumber& \hspace{.5cm} + 2 \gamma \langle \grd_{\bm{\lambda}} {\cal L}( \prm^{(k)}, \bm{\lambda}^{(k)} ) - \grd_{\bm{\lambda}} {\cal L} ( \prm^\star, \bm{\lambda}^\star), \bm{\lambda}^{(k)} - \bm{\lambda}^\star \rangle \\
\nonumber& \hspace{.5cm} + \gamma^2 ( 1 + c_0 ) \| \grd_{\bm{\lambda}} {\cal L}( \prm^{(k)}, \bm{\lambda}^{(k)} ) - \grd_{\bm{\lambda}} {\cal L} ( \prm^\star, \bm{\lambda}^\star) \|^2  \\
& \hspace{.5cm}+ \Big( \frac{2 \gamma}{c_1} + \gamma^2  + \frac{\gamma^2}{c_0 } \Big) \| {\bm e}_{\bm{\lambda}}^{(k)} \|^2.
\end{align}
Summing up the two inequalities gives:
\begin{align}
\nonumber& \| {\bm z}^{(k+1)} - {\bm z}^\star \|^2 \\
\nonumber& \leq (1 + 2 c_1 \gamma ) \!~\| {\bm z}^{(k)} - {\bm z}^\star \|^2 + \Big( \frac{2 \gamma}{c_1} +\gamma^2 +  \frac{\gamma^2}{c_0 } \Big) E_k  \\
\nonumber& \hspace{.25cm} - 2 \gamma \langle \bm{\Phi}( {\bm z}^{(k)} ) - \bm{\Phi}( {\bm z}^\star ),  {\bm z}^{(k)} - {\bm z}^\star \rangle \\
\nonumber& \hspace{.5cm} + \gamma^2 ( 1 + c_0 ) \| \bm{\Phi}( {\bm z}^{(k)} ) - \bm{\Phi}( {\bm z}^\star ) \|^2  \\
\nonumber& \overset{(a)}{\leq}  \Big( 1 + 2 \gamma (c_1  - \upsilon' ) + \gamma^2 ( 1+ c_0) L_{\Phi}^2 \Big) \| {\bm z}^{(k)} - {\bm z}^\star \|^2 \\
& \hspace{.5cm} + \Big( \frac{2 \gamma}{c_1} +\gamma^2 + \frac{\gamma^2}{c_0 } \Big) E_k,
\end{align}
where (a) uses the strong monotonicity and smoothness of the map $\bm{\Phi}$, {\color{blue} cf.~\cite[Lemma 3.4]{koshal2011multiuser}}. 
Setting $c_1 = \upsilon' /2$ yields
\beq
\begin{split}
& \| {\bm z}^{(k+1)} - {\bm z}^\star \|^2 \\
& \leq \Big( 1 - \gamma \upsilon' + \gamma^2 (1+c_0) L_{\Phi}^2 \Big) \| {\bm z}^{(k)} - {\bm z}^\star \|^2 \\
& \hspace{.5cm} +
\Big( \frac{4 \gamma}{\upsilon'} +\gamma^2 + \frac{\gamma^2}{c_0 } \Big) E_k.
\end{split}
\eeq
Observe that we can choose $\gamma$ such that $1 - \gamma \upsilon' + \gamma^2 (1+c_0) L_{\Phi}^2 < 1$.
Moreover, the above inequality implies that $\| {\bm z}^{(k)} - {\bm z}^\star \|^2$ evaluates
to
\beq
\begin{split}
& \| {\bm z}^{(k+1)} - {\bm z}^\star \|^2 \\
& \leq (1 - \gamma \upsilon' + \gamma^2 (1+c_0) L_{\Phi}^2 )^k \| {\bm z}^{(0)} - {\bm z}^\star \|^2 + \\
& \sum_{\ell=1}^k (1 - \gamma \upsilon' + \gamma^2 (1+c_0) L_{\Phi}^2 )^{k-\ell} \Big( \frac{4 \gamma}{\upsilon'} + \gamma^2+ \frac{\gamma^2}{c_0 } \Big) E_\ell.
\end{split}
\eeq
If $E_k \leq \overline{E}$ for all $k$, then ${\bm z}^{(k)}$ converges to a neighborhood of ${\bm z}^\star$ of radius 
\beq
\limsup_{k \rightarrow \infty} \| {\bm z}^{(k)} - {\bm z}^\star \|^2 \leq 
\frac{ \frac{4 \gamma}{\upsilon'} + \gamma^2 +\frac{\gamma^2}{c_0 } }{\gamma \upsilon' - \gamma^2 ( 1+ c_0 ) L_{\Phi}^2} \overline{E}.
\eeq
Setting $c_0 =1$ concludes the proof.
\subsection{Proof of Lemma \ref{lemma:perturbedgradients2}}\label{app:perturbedgradientsinfreq}
Comparing the equations in \eqref{eq:perturbgradients2} with \eqref{eq:prm_update_jam} and \eqref{eq:dl_update_jam}, we identify that:
\begin{equation}
\begin{aligned}
    \big[\boldsymbol{e}^{(k)}_{\prm}\big]_j=\frac{1}{N}\sum_{t=1}^T\lambda_t^{(k)}\big(&\nabla_{\prm}g_t\big({\textstyle\frac{1}{N}\sum_{i=1}^N\widehat{\prm}{}_i^{(k)}}\big)\\
    &-\nabla_{\prm}g_t\big({\textstyle\frac{1}{N}\sum_{i=1}^N\prm_i^{(k)}}\big)\big),
    \end{aligned}
\end{equation}
\begin{equation}
\big[\boldsymbol{e}^{(k)}_{\boldsymbol{\lambda}}\big]_t=g_t\big({\textstyle\frac{1}{N}\sum_{i=1}^N\widehat{\prm}{}_i^{(k)}}\big)-g_t\big(\textstyle{\frac{1}{N}\sum_{i=1}^N\prm_i^{(k)}}\big),
\end{equation}
where $\big[\boldsymbol{e}^{(k)}_{\prm}\big]_j$ denotes the $j$th block of $\boldsymbol{e}^{(k)}_{\prm}$. Using Assumption~\ref{ass:bdd}, we immediately see that:
\begin{equation}
\begin{aligned}
    \| \big[\boldsymbol{e}^{(k)}_{\prm}\big]_j\|&\leq\frac{\overline{\lambda}L T}{N}\|\frac{1}{N}\sum_{i=1}^N(\prm_i^{(k)}-\widehat{\prm}{}_i^{(k)})\|\\
    &\leq \frac{\overline{\lambda}L T}{N^2}\sum_{i=1}^N\|\prm_i^{(k)}-\widehat{\prm}{}_i^{(k)}\|,
    \end{aligned}
\end{equation}
which then implies \eqref{eq:errorthetagrad2}. Assumption~\ref{ass:bdd} implies that $g_t$ is $B$-Lipschitz continuous, therefore
\begin{equation}
    \begin{aligned}
    |\big[\boldsymbol{e}^{(k)}_{\boldsymbol{\lambda}}\big]_t|&\leq B\|\frac{1}{N}\sum_{i=1}^N(\prm_i^{(k)}-\widehat{\prm}{}_i^{(k)})\|\\
    &\leq \frac{B}{N}\sum_{i=1}^N\|\prm_i^{(k)}-\widehat{\prm}{}_i^{(k)}\|,
    \end{aligned}
\end{equation}
which implies \eqref{eq:errorlambdagrad2}.

\subsection{Proof of Lemma \ref{lemma:perturbationupperbound}}\label{app:perturbationbound}
Observe that the gradient perturbation in both dual and primal variables  is upper bounded by some constant times $\sum_{i=1}^N\|\prm_i^{(k)}-\widehat{\prm}{}_i^{(k)}\|$ in \eqref{eq:errorgrad2}. Thus, we would like to upper bound this term. Let $\mathcal{H}_i^{(k)}$ be the set of $(1-\alpha_2)m$ trustworthy parameters of agent $i$ out of the last $m$ parameters at iteration $k$, i.e., $(1-\alpha_2)m$ trustworthy parameters from set $\{\bm r_i^{(k-\ell)}\}_{\ell=0}^{m-1}$. Note that if a parameter is trustworthy, then $\bm r_i^{(k-\ell)}=\prm_i^{(k-\ell)}$. Hence we define the mean of the iterates in set $\mathcal{H}_i^{(k)}$ as:
\begin{equation}
    \overline{\prm}{}_i^{(k)}\eqdef\frac{1}{(1-\alpha_2)m}\sum_{\prm_i^{(k-\ell)}\in \mathcal{H}_i^{(k)}}\prm{}_i^{(k-\ell)}.
\end{equation}
Using triangular inequality, we can write:
\begin{equation}\label{eq:star}
\begin{aligned}
    \|\prm_i^{(k)}-\widehat{\prm}{}_i^{(k)}\|&=\|\prm_i^{(k)}-\overline{\prm}{}_i^{(k)}+\overline{\prm}{}_i^{(k)}-\widehat{\prm}{}_i^{(k)}\|\\
    &\leq \|\prm_i^{(k)}-\overline{\prm}{}_i^{(k)}\|+\|\overline{\prm}{}_i^{(k)}-\widehat{\prm}{}_i^{(k)}\|.
    \end{aligned}
\end{equation}
Let $\widehat{\prm}{}_i^{(k)}$ be the estimated mean using median-based estimator. Using norm equivalence:
\begin{equation}
\begin{aligned}
    \underset{\prm{}_i^{(k-\ell)}\in \mathcal{H}_i^{(k)}}{\max} \|\prm_i^{(k-\ell)}-\overline{\prm}{}_i^{(k)}\|_\infty&\leq \underset{\prm{}_i^{(k-\ell)}\in \mathcal{H}_i^{(k)}}{\max} \|\prm_i^{(k-\ell)}-\overline{\prm}{}_i^{(k)}\|\\
    &\leq \underset{0\leq \ell\leq m-1}{\max} \|\prm_i^{(k-\ell)}-\overline{\prm}{}_i^{(k)}\|.
    \end{aligned}
\end{equation}
Thus, under Assumption~\ref{ass:attackedfraction},  Proposition~\ref{prop:median} suggests:
\begin{equation}\label{eq:estimationerrorbound}
    \|\overline{\prm}{}_i^{(k)}-\widehat{\prm}{}_i^{(k)}\|\leq C_\alpha \underset{0\leq\ell\leq m-1}{\max} \|\prm_i^{(k-\ell)}-\overline{\prm}{}_i^{(k)}\|,
\end{equation}
where $C_\alpha=\frac{2\alpha_2}{1-\alpha_2}\left(1+\sqrt{\frac{(1-\alpha_2)^2}{1-2\alpha_2}}\right)\sqrt{d}$. 

Let $\ell^\star=\underset{0\leq \ell\leq m-1}{\argmax} \|\prm_i^{(k-\ell)}-\overline{\prm}{}_i^{(k)}\|$. Then:
\begin{equation}\label{estimationerrorbound2}
    \begin{aligned}
    &\underset{0\leq \ell\leq m-1}{\max} \|\prm_i^{(k-\ell)}-\overline{\prm}{}_i^{(k)}\|=\|\prm_i^{(k-\ell^\star)}-\overline{\prm}{}_i^{(k)}\|\\
    &\begin{multlined}
        =\|\prm_i^{(k-\ell^\star)}-\prm_i^{(k-\ell^\star+1)}+\prm_i^{(k-\ell^\star+1)}-\dotsc\\
    -\prm_i^{(k-1)}+\prm_i^{(k+1)}-\prm_i^{(k)}+\prm_i^{(k)}-\overline{\prm}{}_i^{(k)}\|
    \end{multlined}\\
    &\leq \|\prm_i^{(k)}-\overline{\prm}{}_i^{(k)}\|+\sum_{j=k-\ell^\star}^{k-1}\|\prm_i^{(j)}-\prm_i^{(j+1)}\|\\
    &\leq \|\prm_i^{(k)}-\overline{\prm}{}_i^{(k)}\|+\sum_{j=k-\ell^\star}^{k-1}\|\gamma [\boldsymbol{\widehat{g}}^{(j)}_{\prm}]_i\|\\
    &\leq \|\prm_i^{(k)}-\overline{\prm}{}_i^{(k)}\|+\gamma\sum_{j=k-m+1}^{k-1}\| \boldsymbol{[\widehat{g}}^{(j)}_{\prm}]_i\|,
    \end{aligned}
\end{equation}
where $\boldsymbol{[\widehat{g}}^{(j)}_{\prm}]_i$ denotes the $i$th block of $\boldsymbol{\widehat{g}}^{(j)}_{\prm}$. Using equations \eqref{eq:estimationerrorbound} and \eqref{estimationerrorbound2}, we can rewrite \eqref{eq:star}:
\beq
\begin{aligned}\label{eq:doublestar}
    \|\prm_i^{(k)}-\widehat{\prm}{}_i^{(k)}\|
    \leq& (1+C_\alpha)\|\prm_i^{(k)}-\overline{\prm}{}_i^{(k)}\|\\&+\gamma C_\alpha\sum_{j=k-m+1}^{k-1}\| \boldsymbol{[\widehat{g}}^{(j)}_{\prm}]_i\|.
\end{aligned}
\eeq
Next step is to bound the $\|\prm_i^{(k)}-\overline{\prm}{}_i^{(k)}\|$ term:
\begin{equation}\label{eq:secondtermbound}
\begin{aligned}
&\|\prm_i^{(k)}-\overline{\prm}{}_i^{(k)}\|=\|\prm_i^{(k)}-\frac{1}{(1-\alpha_2)m}\hspace{-0.1cm}\sum_{\prm{}_i^{(k-\ell)}\in \mathcal{H}_i^{(k)}}\hspace{-.2cm} \prm_i^{(k-\ell)}\|\\
&\leq \frac{1}{(1-\alpha_2)m} \sum_{\prm{}_i^{(k-\ell)}\in \mathcal{H}_i^{(k)}} \|\prm_i^{(k)}-\prm_i^{(k-\ell)}\|\\
&\leq \frac{1}{(1-\alpha_2)m} \sum_{\ell=0}^{m-1} \|\prm_i^{(k)}-\prm_i^{(k-\ell)}\|\\
&\begin{multlined}=\frac{1}{(1-\alpha_2)m} \sum_{\ell=0}^{m-1} \|\prm_i^{(k)}-\prm_i^{(k-1)}+\prm_i^{(k-1)}-\dotsc\\ -\prm_i^{(k-\ell+1)}+\prm_i^{(k-\ell+1)}-\prm_i^{(k-\ell)}\|\end{multlined}\\
&\leq \frac{1}{(1-\alpha_2)m} \sum_{\ell=0}^{m-1}\sum_{j=k-\ell}^{k-1}\|\prm_i^{(j+1)}-\prm_i^{(j)}\|\\
&\leq\frac{m}{(1-\alpha_2)m} \sum_{j=k-m+1}^{k-1}\|\prm_i^{(j+1)}-\prm_i^{(j)}\|\\
&\leq \frac{1}{1-\alpha_2}\gamma \sum_{j=k-m+1}^{k-1}\| \boldsymbol{[\widehat{g}}^{(j)}_{\prm}]_i\|. 
\end{aligned}
\end{equation}
Plugging \eqref{eq:secondtermbound} into \eqref{eq:doublestar}:
\begin{equation}\label{eq:oneagenterrorbound}
\|\prm_i^{(k)}-\widehat{\prm}{}_i^{(k)}\|
    \leq\left(\frac{1+C_\alpha}{1-\alpha_2}+C_\alpha\right)\gamma \sum_{j=k-m+1}^{k-1}\| \boldsymbol{[\widehat{g}}^{(j)}_{\prm}]_i\|.
\end{equation}
For brevity of notation, let $\left(\frac{1+C_\alpha}{1-\alpha_2}+C_\alpha\right)=\overline{C}_\alpha$ and let $\nabla_{\prm}\mathcal{L}_\upsilon^{(k)}:=\nabla_{\prm}\mathcal{L}_\upsilon(\prm^{(k)};\boldsymbol{\lambda}^{(k)})$. Summing up \eqref{eq:oneagenterrorbound} for all agents and using norm equivalence:
\begin{equation}\label{eq:tripplestar}
\begin{aligned}
\sum_{i=1}^N\|\prm_i^{(k)}-\widehat{\prm}{}_i^{(k)}\|
    &\leq\overline{C}_\alpha\sqrt{N}\gamma \sum_{j=k-m+1}^{k-1}\| \boldsymbol{\widehat{g}}^{(j)}_{\prm}\|\\
    &\hspace{-.75cm}\overset{\eqref{eq:thetaperturb2}}{=}\overline{C}_\alpha\sqrt{N}\gamma \sum_{j=k-m+1}^{k-1}\|\nabla_{\prm}\mathcal{L}_\upsilon^{(j)}+\boldsymbol{e}_{\prm}^{(j)}\|\\
    &\hspace{-.75cm}\leq \overline{C}_\alpha\sqrt{N}\gamma \sum_{j=k-m+1}^{k-1}\|\nabla_{\prm}\mathcal{L}_\upsilon^{(j)}\|+\|\boldsymbol{e}_{\prm}^{(j)}\|.
\end{aligned}
\end{equation}
Using \eqref{eq:errorthetagrad2}:
\begin{equation}\label{eq:errorboundoptimalitygap}
    \begin{aligned}
     \|\boldsymbol{e}_{\prm}^{(j)}\|&\leq \frac{\overline{\lambda}LT}{N}\sum_{i=1}^N\|\prm_i^{(j)}-\widehat{\prm}{}_i^{(j)}\|\\
     &=\frac{\overline{\lambda}LT}{N}\sum_{i=1}^N\|\prm_i^{(j)}-\prm_i^\star+\prm_i^\star-\widehat{\prm}{}_i^{(j)}\|\\
     &\leq \frac{\overline{\lambda}LT}{N}\sum_{i=1}^N\|\prm_i^{(j)}-\prm_i^\star\|+\|\widehat{\prm}{}_i^{(j)}-\prm_i^\star\|\\
     &\overset{(*)}{\leq} (1+\sqrt{d})\frac{\overline{\lambda}LT}{N}\sum_{i=1}^N\underset{0\leq \ell_i\leq m-1}{\max}\|\prm_i^{(j-\ell_i)}-\prm_i^\star\|\\
     &\leq (1+\sqrt{d}) \overline{\lambda}LT \underset{\substack{i\\0\leq\ell_i\leq m-1}}{\max}\|\prm_i^{(j-\ell_i)}-\prm_i^\star\|\\
     &\leq  (1+\sqrt{d}) \overline{\lambda}LT \underset{0\leq \ell \leq m-1}{\max}\|\prmdl^{(j-\ell)}-\prmdl^\star\|,
    \end{aligned}
\end{equation}
where $(*)$ is obtained by:
\begin{equation}
\begin{aligned}
     \|\widehat{\prm}{}_i^{(j)}-\prm_i^\star\|&\leq \sqrt{d} \|\widehat{\prm}{}_i^{(j)}-\prm_i^\star\|_\infty\\&\leq \sqrt{d}\underset{0\leq\ell\leq m-1}{\max}\|\prm_i^{(j-\ell)}-\prm_i^\star\|_\infty\\
     &\leq \sqrt{d}\underset{0\leq\ell\leq m-1}{\max}\|\prm_i^{(j-\ell)}-\prm_i^\star\|,
\end{aligned}
\end{equation}
and $\|\prm_i^{(j)}-\prm_i^\star\|\leq\underset{0\leq\ell\leq m-1}{\max} \|\prm_i^{(j-\ell)}-\prm_i^\star\|$.

Furthermore, by using the $L_\Phi$-Lipschitz property of $\boldsymbol{\overline{\Phi}}(\prmdl)$:
\begin{equation}\label{eq:gradientlipschitzbound}
\begin{aligned}
    \|\nabla_{\prm}\mathcal{L}_\upsilon^{(j)}\|\leq\|\boldsymbol{\overline{\Phi}}(\prmdl^{(j)})\|&=\|\boldsymbol{\overline{\Phi}}(\prmdl^{(j)})-\boldsymbol{\overline{\Phi}}(\prmdl^{\star})\|\\
    &\leq \frac{1}{L_\Phi}\|\prmdl^{(j)}-\boldsymbol{z^\star}\|\\
    &\leq \frac{1}{L_\Phi}\underset{0\leq \ell \leq m-1}{\max}\|\prmdl^{(j-\ell)}-\boldsymbol{z^\star}\|
    \end{aligned}
\end{equation}
We can rewrite \eqref{eq:tripplestar} using equations \eqref{eq:errorboundoptimalitygap} and \eqref{eq:gradientlipschitzbound}:
\begin{equation}
\begin{aligned}
    &\sum_{i=1}^N\|\prm_i^{(k)}-\widehat{\prm}{}_i^{(k)}\|
    \leq\overline{C}_\alpha\sqrt{N}\gamma \sum_{j=k-m+1}^{k-1}\|\nabla_{\prm}\mathcal{L}_\upsilon^{(j)}\|+\|\boldsymbol{e}_{\prm}^{(j)}\|\\
    &
    \leq \overline{C}_\alpha\sqrt{N} C_0 \gamma \sum_{j=k-m+1}^{k-1}\underset{0\leq \ell \leq m-1}{\max}\|\prmdl^{(j-\ell)}-\boldsymbol{z^\star}\|\\
    &\leq \overline{C}_\alpha\sqrt{N}C_0 (m-1) \gamma \underset{1\leq \ell \leq 2(m-1)}{\max}\|\prmdl^{(k-\ell)}-\boldsymbol{z^\star}\|,
    \end{aligned}
\end{equation}
where $C_0=\frac{1}{L_\Phi}+(1+\sqrt{d})\overline{\lambda}LT$. Finally, using \eqref{eq:errorgrad2} and letting $C_1=\left(\frac{\overline{\lambda}LT}{N}\right)^2+\left(\frac{BT}{N}\right)^2$:
\begin{equation}\begin{aligned}
    E_k&=\|\boldsymbol{e}^{(k)}_{\prm}\|^2+\|\boldsymbol{e}^{(k)}_{\boldsymbol{\lambda}}\|^2\leq C_1\left(\sum_{i=1}^N\|\prm_i^{(k)}-\widehat{\prm}{}_i^{(k)}\|\right)^2\\
    &\leq C_1 (\overline{C}_\alpha\sqrt{N}C_0 (m-1))^2\gamma^2 \hspace{-.1cm}\underset{1\leq \ell \leq 2(m-1)}{\max}\hspace{-.1cm}\|\prmdl^{(k-\ell)}-\boldsymbol{z^\star}\|^2\\
    &\leq C_1 (\overline{C}_\alpha\sqrt{N}C_0 (m-1))^2\gamma^2 \hspace{-.1cm}\underset{0\leq \ell \leq 2(m-1)}{\max}\hspace{-.1cm}\|\prmdl^{(k-\ell)}-\boldsymbol{z^\star}\|^2\\
    &\leq \overline{C}\gamma^2\underset{0\leq \ell \leq 2(m-1)}{\max}\|\prmdl^{(k-\ell)}-\boldsymbol{z^\star}\|^2,
    \end{aligned}
\end{equation}
where $\overline{C}=\left(\frac{T^2(\overline{\lambda}^2L^2+B^2)}{N}\right)\times\left(\frac{1}{L_\Phi}+(1+\sqrt{d})\overline{\lambda}LT\right)^2\times$ $\left(\frac{1+C_\alpha}{1-\alpha_2}+C_\alpha\right)^2\times(m-1)^2$.
\subsection{Proof of Theorem \ref{thm:convergence}}\label{app:infrequentthm}
Based on Lemma~\ref{lemma:perturbedgradients2}, our idea is to perform a perturbation analysis on the PDA algorithm. The first part of the proof is analogous to that of Theorem~\ref{thm:main}, and then upper bounding $E_k$ by Lemma~\ref{lemma:perturbationupperbound}. This yields:
\begin{equation}\label{eq:thm2pf1}
\begin{aligned}
    &\|\prmdl^{(k+1)}-\prmdl^\star\|\leq (1-\gamma \upsilon+2\gamma^2L_\Phi^2)\|\prmdl^{(k)}-\prmdl^\star\|^2\\
    &+\left(\frac{4\gamma}{\upsilon}+2\gamma^2\right)\gamma^2\overline{C}\underset{0\leq \ell\leq2(m-1)}{\max}\|\prmdl^{(k-\ell)}-\prmdl^\star\|^2.
\end{aligned}
\end{equation}

For the second part of the proof, we use the following Lemma:
\begin{lemma}\label{lem:johansson}
\cite[Lemma 3]{mjohansson} Let $\{V(t)\}$ be a sequence of real numbers satisfying
\begin{equation*}
    V(t+1)\leq pV(t)+q\underset{t-\tau(t)\leq s\leq t}{\max}V(s)+r,\quad t\in \mathbb{N}_0,
\end{equation*}
for some nonnegative constants $p$,$q$, and $r$. If $p+q<1$ and
\begin{equation*}
    0\leq\tau(t)\leq\tau_{\max},\quad t\in \mathbb{N}_0,
\end{equation*}
then
\begin{equation*}
    V(t)\leq \rho^tV(0)+\epsilon,\quad t\in \mathbb{N}_0,
\end{equation*}
where $\rho=(p+q)^{\frac{1}{1+\tau_{\max}}}$ and $\epsilon=r/(1-p-q)$.
\end{lemma}
We apply Lemma \ref{lem:johansson} on \eqref{eq:thm2pf1} for  $t=k\geq2(m-1)$, $V(t)=V(k)=\|\prmdl^{(k)}-\prmdl^\star\|^2$, $p=1-\gamma \upsilon+2\gamma^2L_\Phi^2$, $q=\left(\frac{4\gamma}{\upsilon}+2\gamma^2\right)\gamma^2\overline{C}$, $r=0$, and $\tau_{\max}=2(m-1)$ to get:
\begin{equation}\label{eq:geometricconvpf}
    V(k)\leq\rho^{k-2(m-1)}V(2(m-1)),\quad k\geq2(m-1),
\end{equation}
where $\rho=(1-\gamma \upsilon+2\gamma^2L_\Phi^2+\frac{4\overline{C}\gamma^3}{\upsilon}+2\overline{C}\gamma^4)^\frac{1}{1+2(m-1)}$. The condition $p+q<1$ is met when:
\begin{equation*}
    f(\gamma)=\upsilon-2\gamma L_\Phi^2-\frac{4\overline{C}\gamma^2}{\upsilon}-2\overline{C}\gamma^3>0
\end{equation*}
Observe that $f(\gamma)$ is a continuous function in $\gamma$, and $f(0)=\upsilon>0$. Hence, there exists a small $\gamma>0$ such that $f(\gamma)>0$, which satisfies the required condition. Taking the limit as $k$ goes to infinity in \eqref{eq:geometricconvpf}:
\begin{equation}
    \underset{k\rightarrow\infty}{\lim}V(k)\leq \underset{k\rightarrow\infty}{\lim}\rho^{k-2(m-1)}V(2(m-1))=0,
\end{equation}
since $\rho<1$. Finally, since $V(k)\geq 0$, we conclude that $\underset{k\rightarrow\infty}{\lim}V(k)=\underset{k\rightarrow\infty}{\lim}\|\prmdl^{(k)}-\prmdl^\star\|^2=0$.
\bt{\subsection{Proof of Bounded Dual Variables}\label{app:boundedgradient}
The proof is the same for statements in both Lemma~\ref{lem:perturbedgradients1} and Lemma~\ref{lemma:perturbedgradients2}.
The update rule given by \eqref{eq:dual_pb}:
		\beq\label{eq:pfupdatedualstatic}\lambda_t^{(k+1)} = \big[ \lambda_t^{(k)} + \gamma \big( \overline{g}_t ((1-\alpha_1) \widehat{\prm}_{\cal H}^{(k)} ) - \upsilon \lambda_t^{(k)} \big) \big]_+. \eeq
The update rule given by \eqref{eq:dl_update_jam}:
\beq\label{eq:pfupdatedualdynamic}
\lambda_t^{(k+1)} = \big[ \lambda_t^{(k)} + \gamma \big( g_t ( {\textstyle \widehat{\prm}^{(k)} ) - \upsilon \lambda_t^{(k)} \big) \big]_+}. 
\eeq
		Let $\overline{g}_t(\cdot)\leq M$ and $0\leq\lambda_t^{(k)}\leq \frac{M}{\upsilon}$. We can upper bound both \eqref{eq:pfupdatedualstatic} and \eqref{eq:pfupdatedualdynamic} as:
		\beq\lambda_t^{(k+1)}\leq (1-\gamma\upsilon)\lambda_t^{(k)}+\gamma M\leq \frac{M}{\upsilon}.\eeq
		Hence, one can set $\lambda_t^{(0)}\leq \frac{M}{\upsilon}$ to guarantee assumption. Anyhow
		if $\frac{M}{\upsilon}\leq\lambda_t^{(k)}$:
		\beq\lambda_t^{(k+1)}\leq (1-\gamma\upsilon)\lambda_t^{(k)}+\gamma M\leq \lambda_t^{(k)}.\eeq
		Thus, if $\lambda_t^{(0)}\leq \frac{M}{\upsilon}$ then $\overline{\lambda}=\frac{M}{\upsilon}$. If $\lambda_t^{(0)}\geq \frac{M}{\upsilon}$ then $\overline{\lambda}=\lambda_t^{(0)}$. This guarantees the assumption.}

\end{document}